\documentclass[aip,
 amsmath,
 amssymb,
 reprint,
 floatfix
 ]{revtex4-1}

\draft 

\usepackage{lipsum}
\usepackage{amsfonts}
\usepackage{graphicx}
\usepackage{epstopdf}
\usepackage{mathdots}
\usepackage{nicefrac}
\usepackage{enumitem}

\usepackage{algorithm, algpseudocode}
\algrenewcommand\algorithmicrequire{\textbf{Input:}}
\algrenewcommand\algorithmicensure{\textbf{Output:}}

\usepackage{amsthm}
\usepackage{bm}
\usepackage{bbm}
\usepackage{color}

\newtheorem{theorem}{Theorem}[section]
\newtheorem{conjecture}[theorem]{Conjecture}
\newtheorem{corollary}[theorem]{Corollary}
\newtheorem{lemma}[theorem]{Lemma}
\newtheorem{remark}[theorem]{Remark}
\newtheorem{hypotheses}[theorem]{Hypotheses}

\newcommand{\sympmap}{\bm{F}} 
\newcommand{\obs}{\tilde{\bm{h}}}
\newcommand{\dobs}{\tilde{\bm{g}}}
\newcommand{\obsT}{\bm{h}}
\newcommand{\dobsT}{\bm{g}}
\newcommand{\BA}{\mathcal{B}} 
\newcommand{\WBA}{\mathcal{W}\mathcal{B}} 
\newcommand{\Kback}{\bar{K}} 
\newcommand{\Whalf}{\tilde{W}_K}
\newcommand{\refpoly}{q_{\alpha, n}}

\newcommand{\abf}{\bm{a}} 
\newcommand{\Bbf}{\bm{B}} 
\newcommand{\filter}{\bm{c}} 
\newcommand{\halffilter}{\tilde{\bm{c}}} 
\newcommand{\dbf}{\bm{d}} 
\newcommand{\ebf}{\bm{e}} 
\newcommand{\nbf}{\bm{n}} 
\newcommand{\ubf}{\bm{u}} 
\newcommand{\xbf}{\bm{x}} 
\newcommand{\zbf}{\bm{z}} 
\newcommand{\thetabf}{\bm{\theta}}
\newcommand{\omegabf}{\bm{\omega}}
\newcommand{\ksm}{k_{\mathrm{sm}}} 

\newcommand{\xibf}{\bm{\xi}} 

\newcommand{\Nbb}{\mathbb{N}} 
\newcommand{\Rbb}{\mathbb{R}} 
\newcommand{\Tbb}{\mathbb{T}} 
\newcommand{\Zbb}{\mathbb{Z}} 

\newcommand{\Fcal}{\mathcal{F}} 
\newcommand{\Ocal}{\mathcal{O}} 

\newcommand{\abs}[1]{\left | #1 \right |}
\newcommand{\norm}[1]{\left \lVert #1 \right \rVert}

\newcommand{\dif}{\hspace{1pt} \mathrm{d}}
\newcommand{\ceil}[1]{\left\lceil #1 \right\rceil}
\newcommand{\floor}[1]{\left\lfloor #1 \right\rfloor}

\DeclareMathOperator{\Tr}{\mathrm{Tr}}

\begin{document}


\title{Finding Birkhoff Averages via Adaptive Filtering} 



\author{M. Ruth}
\email[]{mer335@cornell.edu}
\affiliation{Center for Applied Mathematics, Cornell University, Ithaca, NY}
\author{D. Bindel}
\affiliation{Department of Computer Science, Cornell University, Ithaca, NY }


\date{\today}

\begin{abstract}
In many applications, one is interested in classifying trajectories of Hamiltonian systems as invariant tori, islands, or chaos. The convergence rate of ergodic Birkhoff averages can be used to categorize these regions, but many iterations of the return map are needed to implement this directly. Recently, it has been shown that a weighted Birkhoff average can be used to accelerate the convergence, resulting in a useful method for categorizing trajectories.

In this paper, we show how a modified version the reduced rank extrapolation method (named Birkhoff RRE) can also be used to find optimal weights for the weighted average with a single linear least-squares solve.
Using these, we classify trajectories with fewer iterations of the map than the standard weighted Birkhoff average. 
Furthermore, for the islands and invariant circles, a subsequent eigenvalue problem gives the number of islands and the rotation number.
Using these numbers, we find Fourier parameterizations of invariant circles and islands. 
We show examples of Birkhoff RRE on the standard map and on magnetic field line dynamics.
\end{abstract}


\maketitle 

\section{Introduction}
Invariant tori are ubiquitous structures in symplectic maps and Hamiltonian dynamics. 
Examples of invariant tori include periodic orbits of the pendulum, invariant circles in the standard map, halo orbits in astrodynamics \cite{kolemen2012}, and nested flux surfaces in magnetic confinement devices\cite{paul2022}.
Such orbits are known to be stable to perturbation due to the KAM theorem \cite{LLave_2004}.
However, numerically identifying orbits from trajectories is often challenging, due to both the difficulties of high-dimensional geometry and the problems of small denominators.

One standard numerical method to find invariant tori is the parameterization method \cite{haro2016}.
This method is based on the conjugacy relation defining invariant tori, and can be accelerated using the fast Fourier transform.
This allows for highly accurate computations of invariant tori, which can be proven to exist by a numerical variant of the KAM theorem \cite{Figueras2017}. 
However, one of the main drawbacks of the parameterization method is that it needs an initial guess.
In the case of 1D and 1.5D Hamiltonian systems, a manual initial guess can found relatively easily with a Poincar\'e plot of a trajectory.
Once one solution is found, continuation\cite{kolemen2012} can be used to find more solutions.
However, in the cases of higher dimensional systems or islands, initial guesses are significantly harder to find.
The initial guess issue is faced by other methods relying on iterations on torus parameterizations, such as the flux minimizing surfaces\cite{dewar1992,dewar1994}.

Additionally, initial guesses of the rotation number are also difficult.
Simple methods for finding the rotation are available for invariant circles where there is a natural point to wind about.
In such cases, one can find the rotation number via classical limits\cite{guckenheimer1997} or more accurately using the weighted Birkhoff average\cite{Das2018,sander2020,sander2023} or Richardson-like extrapolation algorithms\cite{Luque2014,Villanueva2022}. 
Once the rotation number is found, then it is straightforward to find a parameterization of the orbit \cite{blessing2023}.
Without good guesses at the winding structure, the typical solution is to look for peaks in the Fourier spectrum of the signal or to use another frequency-based method\cite{laskar1999}. 
Unfortunately, these peaks are only as accurate as the discretization resolution of the spectrum.
This issue is again made more complicated in higher dimensions, where winding is a less well-defined concept.

An additional issue that we will discuss is that of orbit classification.
Before an invariant torus can be fit to a trajectory, we must first be confident that trajectory is, in fact, a torus.
Classical methods for this typically rely on the Lyapunov exponent, but this can be quite slow to converge.
More recently, the rapid convergence of the weighted Birkhoff average\cite{Das2016,Das2017,Das2018,Das2019,sander2020,sander2023} (WBA) has been shown to be capable distinguishing chaotic from non-chaotic trajectories.

Symplectic maps are often computed by evolving Hamiltonian dynamics by numerical integration, and this can present its own problems.
One problem is noise: reliable methods must be robust to potentially non-symplectic errors in the symplectic map.
This is particularly relevant when symplectic integrators are not available for a given application.
We note that for the parameterization method, this is alleviated by an overdetermined formulation of the method\cite{Baresi2018}.
A second issue is the cost of evaluating the map.
Both the parameterization method using the Fast Fourier Transform and WBA can be performed in nearly linear time in the number of samples, so the dominant cost is typically dominated by evolving Hamiltonian trajectories.
When classifying many trajectories, any reduction in the number of evaluations is very useful.
Structured symplectic interpolants such as the HenonNet\cite{Burby2021}, SympNet\cite{Jin2020}, or a Gaussian Process approach\cite{Rath2021} can alleviate both issues of numerical integration.
However, symplectic interpolants are necessarily nonlinear, and may require more data than simply finding the invariant circle. 
Linear interpolants can similarly be used, but require high accuracy, and the evaluation time of constructing the interpolant may still represent the dominant cost of any algorithm.

In this paper, we classify and parameterize invariant tori from single trajectories efficiently in the number of map evaluations. 
For the classification step, we will rely on a variant of the reduced rank extrapolation method\cite{Sidi2017}, which we will refer to as Birkhoff RRE.
Birkhoff RRE works by attempting to find a filter (or linear model) for measurements on a trajectory using only the time-series information.
Because Birkhoff RRE only depends on a single trajectory, it does not require any initial guesses or continuation of invariant tori.
Additionally, Birkhoff RRE is written as a linear least-squares problem, meaning its implementation is straightforward and there is a residual indicating the fit of the linear model.
We prove that when the trajectory is on an invariant torus, the residual of Birkhoff RRE approaches zero as rapidly as the weighted Birkhoff average, allowing for the same classification ability.
Additionally, we show experimentally that RRE converges significantly faster than WBA on a set of trajectories of the standard map, with a large majority of trajectories being classified to highly accurate residuals below $10^{-11}$ in fewer than $1000$ iterations of the map.

For the parameterization step, we show experimentally that the frequencies filtered by the Birkhoff RRE filter are multiples of the rotation number. 
Using those frequencies, we numerically identify both the number of islands and rotation numbers for trajectories in 2D. 
Once this information is known, we project the signal back onto the corresponding Fourier modes, giving a parameterization of the invariant circle or island.
This process is similar to the filter diagonalization method\cite{Mandelshtam2001}.

We introduce necessary background in Sec.~\ref{sec:background}, including a review of the weighted Birkhoff average.
In Sec.~\ref{sec:method}, we build upon WBA to introduce the Birkhoff RRE algorithm, stating the relevant convergence theorems.
Then, in Sec.~\ref{sec:examples}, we show two examples of the method.
In the first example, we examine the convergence of Birkhoff RRE on the standard map, confirming the predicted rates from Sec.~\ref{sec:method}.
Then, in the second example, we show how the method can be applied to an example on a symplectic map obtained from a toroidal plasma confinement device known as a stellarator, showing how the method can be used in a real-world situation.
Finally, we conclude in Sec.~\ref{sec:conclusions}.

\section{Background}
\label{sec:background}

\begin{figure*}
    \centering
    \includegraphics[width=0.7\textwidth]{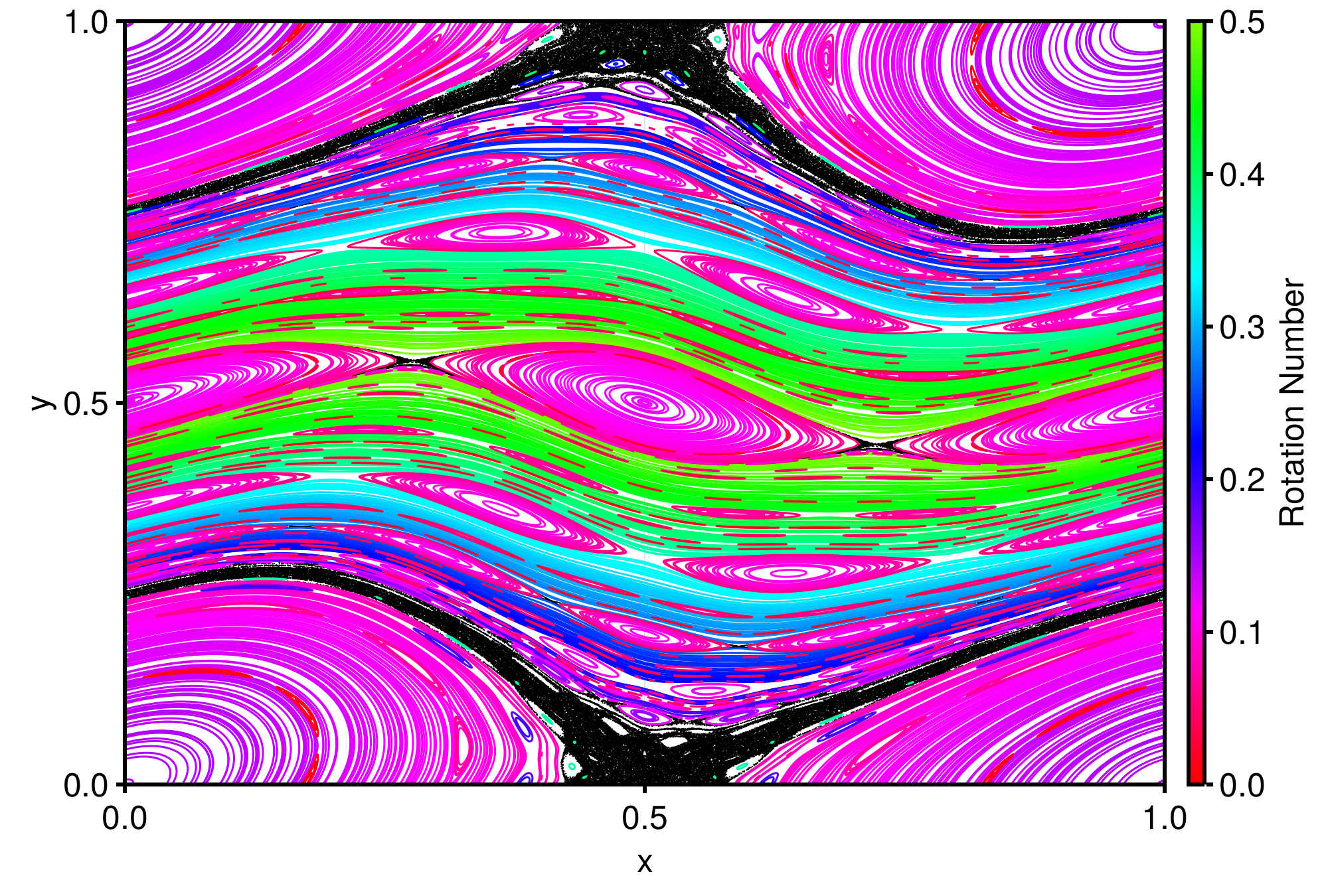}
    \caption{Phase portrait of the standard map \eqref{eq:standard-map}. Parameterizations of the invariant circles and islands are obtained via the methodology in Sec.~\ref{sec:method}, using Algorithm \ref{alg:adaptiveRRE} with $\epsilon=0$, $\gamma=3$, $\delta=10^{-10}$, $K_{\mathrm{init}}=50$, $K_{\max}=600$, and $\Delta K = 50$. Invariant circles and islands are colored according to their rotation numbers, while trajectories classified as chaotic are plotted in black. }
    \label{fig:standard-map}
\end{figure*}

Let $\sympmap : X \to X$ be a map for a discrete-time dynamical system, where $X$ is a $C^M$ manifold. A $d$-dimensional \textit{invariant torus} of $\sympmap$ is a function $\zbf : \Tbb^d \to X$ that satisfies the conjugacy
\begin{equation*}
    \sympmap \circ \zbf - \zbf \circ R_\omega = 0,
\end{equation*}
where $R_\omega : \Tbb^d \to \Tbb^d$ is the rotation map $R_\omega(\thetabf) = \thetabf + \omegabf$ with $\omegabf \in \Tbb^d$. The invariant torus has an associated invariant measure $\mu$ on its image $X_0 \subseteq X$. We assume that the invariant tori are smooth, i.e.~$\zbf \in C^M$ for some positive integer $M$ or $\zbf \in C^\infty$. We refer to $d=1$ invariant tori as \textit{invariant circles}. 

We require the rotation vector $\omegabf$ satisfy the $(c, \nu)$ Diophantine condition 
\begin{equation*} 
    \abs{\omegabf \cdot \nbf - m} \geq \frac{c}{\norm{\nbf}^{\nu}} \qquad \text{for all } \nbf\in\Zbb^d\backslash \{0\}, \ m \in \Zbb.
\end{equation*}
For $d=1$, the Diophantine condition tells us the $\omega \in \Rbb$ is ``sufficiently irrational.''
In higher dimensions, the condition gives a measure of irrational independence, requiring that the rotation numbers also be far from being rational multiples of the others.
We note that the Diophantine condition is satisfied by almost all rotation numbers under the Lebesgue measure \cite{LLave_2004}.
So, almost all invariant tori are Diophantine for nested regions with shear. 

We define an \textit{island} to be a set of $p \geq 1$ tori $\zbf^{(j)}:\Tbb^d \to X$ with $1 \leq j \leq p$ such that for some $\omegabf$ the following conjugacy is satisfied:
\begin{equation}
\label{eq:island-conjugacy}
    \sympmap \circ \zbf^{(j)} = \begin{cases}
        \zbf^{(j+1)}, & 1 \leq j < p,\\
        \zbf^{(1)} \circ R_\omega, & j = p.
    \end{cases}
\end{equation}
A direct result of the above definition is that each $\zbf^{(j)}$ is an invariant torus of the map $\sympmap^p$ (and hence a $p=1$ island is equivalent to an invariant torus). 
Another special case of an island is a periodic orbit, where each $\zbf^{(j)}$ is a constant function.
If a trajectory is on an invariant torus or island, we call it \textit{integrable}.

We note that the definitions of invariant tori and islands above do not require any special properties of the map.
However, they are both commonly found in the case of \textit{symplectic maps}, which preserve some symplectic 2-form under the pushforward of the map.
In the case of symplectic maps, we refer to everything that is not integrable as \textit{chaotic}. 
We note that this is a heuristic definition of chaos, as the theorems herein only guarantee convergence rates of certain methods for integrable trajectories.
The converse (i.e.~chaotic trajectories converge slowly for some definition of chaos) is still an open problem.

In Fig.~\ref{fig:standard-map}, we plot the phase portrait of the Chirikov standard map $F : (x_t, y_t) \mapsto (x_{t+1}, y_{t+1})$  on $X = \Tbb \times \Rbb$ where
\begin{align}
\nonumber
    x_{t+1} &= x_t + y_{t+1} \mod 1,\\
\label{eq:standard-map}
    y_{t+1} &= y_t - \frac{\ksm}{2\pi} \sin (2\pi x_t),
\end{align}
and $\ksm=0.7$. The map has invariant circles (e.g.~nested about $(0,0)$ and the blue-to-green gradient of trajectories in the center), islands (e.g.~the $p=2$ island chain centered at $(0.0,0.5)$ and $(0.5,0.5)$), and chaos in black (e.g.~about $(0.5, 0)$). The invariant circles and islands are both colored by the rotation number $\omega$, found using the methods in Sec.~\ref{sec:method}. We note that $\omega$ has a unique representation in $[0.0,0.5]$ due both to the fact that $\omega \in \Tbb$ and the freedom to take $\theta \to -\theta$ in the parameterization.

To categorize orbits, we consider the problem of finding ergodic averages. Let $\obs : X \to \Rbb^D$ be an observable function on our state space. The \textit{Birkhoff average} of $\obs$ is defined as the limit of finite time averages:
\begin{equation*}
    \BA[\obs](\xbf) = \lim_{\Kback \to \infty} \BA_{\Kback}[\obs](\xbf)
\end{equation*}
where
\begin{equation*}
    \BA_{\Kback}[\obs] = \frac{1}{\Kback} \sum_{k=0}^{\Kback-1} (\obs \circ \sympmap^k)(\xbf).
\end{equation*}
For an intial point $\xbf$, we let $X_0$ be the ergodic component
\begin{multline*}
    X_0 = \{\xbf' \in X \mid \forall f \in C^b(X), \\ \lim_{T \to \infty} \BA_{T}[f](\xbf) - \BA_{T}[f](\xbf') = 0 \}
\end{multline*}
where $C^b(X)$ is the set of continuous bounded functions on $X$. 
We note that the ergodic component $X_0$ is identical to our previous definition for circles with irrational rotation numbers $\omega$.
Then, for almost all $\xbf$, the Birkhoff average converges to an average over a unique invariant measure $\mu$ on $X_0$\cite{coudene_ergodic_2016} 
\begin{equation*}
    \BA[\obs](\xbf) = \int_{X_0}\obs(\xbf) \dif \mu.
\end{equation*}
Additionally, when $\obs\in C^M$ for $M>\nu+d$, one can show that the partial averages have an error $\abs{\BA[\obs] - \BA_{\Kback}[\obs]} = \Ocal(\Kback^{-1})$ on invariant circles and islands. In contrast, chaotic trajectories are conjectured to have a convergence rate of $\Ocal(\Kback^{-1/2})$\cite{sander2020}, the same convergence as expected from a central limit theorem.

If we compose the observable with an invariant circle, we can define the observable using coordinates on the torus. That is, if we let $\obsT = \obs \circ \zbf$, we define the finite-time Birkhoff average of the function $\obsT$ at a point $\thetabf \in \Tbb^d$ as 
\begin{align*}
    \BA_{\Kback}[\obsT](\thetabf) &= \frac{1}{\Kback} \sum_{k=0}^{\Kback-1}(\obsT \circ R_\omega^k)(\thetabf),\\
    &= \BA_{\Kback}[\obs](\zbf(\thetabf)).
\end{align*}
Assuming $\obsT$ is continuous and $\omegabf$ is irrational and rationally independent (i.e.~$R_\omega$ is ergodic on $\Tbb^d$), the limit of these averages is independent of the initial point $\thetabf$. The average is equal to a spatial average
\begin{equation}
\label{eq:BirkhoffAverage}
    \BA[\obsT] = \lim_{\Kback \to \infty} \BA_{\Kback}[\obsT] = \int_{\Tbb^d} \obsT(\thetabf) \dif \thetabf,
\end{equation}
where we note that the Lebesgue measure is the unique invariant measure under $R_\omega$. 
In this way, we connect averages of time series to averages over invariant tori.

One application of Birkhoff averages is to find Fourier coefficients of $\obsT$. Consider that
\begin{equation*}
    \obsT(\thetabf) = \sum_{\nbf \in \Zbb^d} \obsT_{\nbf} e^{2\pi i \nbf \cdot \thetabf}.
\end{equation*}
Then, the coefficients $\obsT_{\nbf}$ are determined by
\begin{align}
\label{eq:FourierBirkhoff}
    \obsT_{\nbf} &= \int_{\Tbb^d} \obsT(\theta)e^{-2\pi i \nbf \cdot \thetabf} \dif \thetabf \\
\nonumber
    &= \BA\left[\obsT(\star) e^{-2\pi i \nbf \cdot \star}\right]\\
\nonumber
    &= \lim_{\Kback \to \infty} \frac{1}{\Kback} \sum_{k=0}^{\Kback-1}\obsT(k \omegabf) e^{2\pi i k \nbf \cdot \omegabf},
\end{align}
where we use `$\star$' for the argument to be averaged over.
A special case is $\nbf = 0$, where the constant Fourier term aligns with the unweighted Birkhoff average $\obsT_0 = \BA[\obsT]$. 
If $\obs$ is the identity and $X$ is a Euclidean space, the coefficients $\obsT_{\nbf} = \zbf_{\nbf}$ provide a Fourier parameterization of an invariant torus
\begin{equation*}
    \zbf(\thetabf) = \sum_{n\in\Zbb} \zbf_{\nbf} e^{2\pi i \nbf \cdot \thetabf} = \sum_{n\in\Zbb}\BA\left[\zbf(\star) e^{-2\pi i \nbf \cdot \omegabf \star}\right] e^{2\pi i n \theta}.
\end{equation*}

To apply the above process for finding coefficients, one must first obtain the rotation vector. For invariant circles, another application of the ergodic average is to obtain the rotation number $\omega$. For instance, if one has access to a symplectic map $\sympmap:\Tbb\times \Rbb \to \Tbb \times \Rbb$ of the form 
\begin{equation}
\label{eq:winding_map}
    x_{t+1} = x_t + y_{t+1} \mod 1, \qquad y_{t+1} = F_y(x_t, y_t),
\end{equation}
the average of $y$ is the rotation number of irreducible circles (i.e., those that wind around the torus). Additionally, if one has access to a point $\xbf_0$ inside of an invariant circle $\zbf \subset \Rbb^2$, one can often find the rotation number by averaging the winding about that point.

However, the winding process to find rotation numbers has some potential difficulties. 
The most immediate difficulty is with orbits that have complicated shapes (such as crescent- or banana-like orbits), where the average position lies outside of the circle.
Even when a point is known to be within the circle, it must be star-shaped for winding to be successful.
A more difficult situation is when $\zbf$ is not injective, as can occur when considering delay embeddings.

With some modifications, we can modify the above ideas to analyze islands.
Using the fact that islands consist of $p$ periodic circles, one can similarly define Fourier series for the observable on each island $\obsT^{(j)} = \obs\circ \zbf^{(j)}$. 
Then, we see that the subsequences associated with each island in the chain can be written as
\begin{equation}
\label{eq:a-island-single}
    \abf_{j+k p} = \sum_{\nbf\in\Zbb^d} \obsT_{\nbf}^{(j)}e^{2\pi i k \nbf \cdot \omegabf}.
\end{equation}
For each $n$, the coefficients $\obsT_{\nbf}^{(j)}$ can be written as a finite discrete Fourier series 
\begin{equation*}
    \obsT^{(j)}_{\nbf} = \sum_{m = 0}^{p-1} \obsT_{m\nbf} e^{2\pi i m j / p}.
\end{equation*}
Combining this with \eqref{eq:a-island-single}, we find that the full sequence can be written as 
\begin{equation}
\label{eq:a-island}
    \abf_t = \sum_{\nbf\in\Zbb^d} \sum_{m=0}^{p-1} \obsT_{m\nbf} e^{2\pi i t (\nbf \cdot \omegabf + m)/p}
\end{equation}
So, signals associated with islands have two frequency components: the rational frequencies $m/p$ associated with jumping between islands and the irrational frequencies $\nbf \cdot \omegabf/p$ associated with the rotation number.
A consequence is that for the map \eqref{eq:winding_map}, averaging $y_t$ will typically return a rational number $m/p$ associated with the number of islands, rather than $\omega$. 
To find $\omega$ for islands, a two-step process would then be needed, where first the denominator of the average of $y_t$ is to identify islands (as performed in Sander and Meiss\cite{sander2020}), and then use a second average of $\sympmap^p$ is used to determine rotation.
In the case of higher-dimensional tori and islands, there are multiple irrational frequencies, making it increasingly difficult to solve for frequencies using winding. 

Another difficulty is that the $\Ocal(\Kback^{-1})$ convergence rate of the Birkhoff averages on circles and islands is too slow for most applications.
For smooth enough maps and invariant circles, this can be improved via the \textit{weighted Birkhoff average}
\begin{equation}
\label{eq:WBA}
    \WBA_{\Kback}[\obs](\xbf) = \sum_{k = 0}^{\Kback-1} w_{k,\Kback} (\obs \circ \sympmap^k)(\xbf),
\end{equation}
where the coefficients $w_{k,\Kback}$ are sampled from a positive function $w\in C^{\infty}$ compactly supported on $[0,1]$ as
\begin{equation*}
    w_{k,\Kback} = \left(\sum_{j=0}^{\Kback-1} w\left(\frac{j+1}{\Kback+1} \right)\right)^{-1} w\left(\frac{k+1}{\Kback+1} \right).
\end{equation*}

To state the convergence rate theorem of the weighted Birkhoff average, we first summarize our assumptions:
\begin{hypotheses}
\label{setting}
    We assume that the three hypotheses hold:
    \begin{enumerate}[label=H\arabic*.]
        \item (Smooth bump function) Let $w\in C^{\infty}(\Rbb)$ have compact support on $[0,1]$ and $w(x)>0$ for all $x\in(0,1)$.
        \item (Diophantine) Let $\omegabf$ be a rotation vector satisfying a $(c,\nu)$ Diophantine condition.
        \item (Smooth system) Let $\sympmap : X \to X$ be a map where $X$ is a $C^M$ manifold and $\sympmap \in C^M$. Additionally, let $X_0 \subseteq X$ be a set where $\sympmap$ is conjugate to island dynamics with $\zbf^{(j)}\in C^M$ (see \eqref{eq:island-conjugacy}), rotation vector $\omegabf$, invariant measure $\mu$ and period $p\geq1$. Finally, let $\obs:X\to \Rbb^D$ be an observable in $C^M$.
    \end{enumerate}
\end{hypotheses}
Additionally, we note that the simpler case of an invariant torus can be considered by setting $p=1$ in the above setting. The following theorem gives the convergence rate of the weighted Birkhoff average:
\begin{theorem}[Das and Yorke\cite{Das2018} Thm 3.1] 
\label{thm:Das2018}
Let $m>1$ be an integer. Under Hypotheses~\ref{setting}, there is a constant $C_m$ depending on $w$, $\obs$, $m$, $M$, $\nu$, and $p$ but independent of $\xbf \in X_0$ such that
\begin{equation*}
    \abs{\WBA_{\Kback} [\obs](\xbf) - \int_{X_0}\obs \dif \mu} \leq C_m \Kback^{-m},
\end{equation*}
provided the `smoothness' $M$ satisfies
\begin{equation*}
    M > d + m\nu.
\end{equation*}
\end{theorem}
We have modified Thm. \ref{thm:Das2018} for this work by assuming $w\in C^\infty$ and extending the theorem to the case of islands. 
The proof that the theorem works for islands is a simple extension of the invariant torus case, whereby the frequencies in \eqref{eq:a-island} are used instead of those in the original proof.
A less smooth function $w$ could be considered, but we have not found this to be useful in practice.

When $M=\infty$ for Thm.~\ref{thm:Das2018}, one can obtain constants $C_m$ for any $m \in \Nbb$. 
That is, the weighted Birkhoff average converges as $\Ocal(\Kback^{-m})$ for all $m$. 
This is useful if one wants to quickly compute Birkhoff averages. 
In contrast, this result does not hold for chaotic trajectories, so the convergence rate of the weighted Birkhoff average can be used to classify trajectories \cite{sander2020}. 

\begin{figure*}
    \centering
    \includegraphics[width=0.8\textwidth]{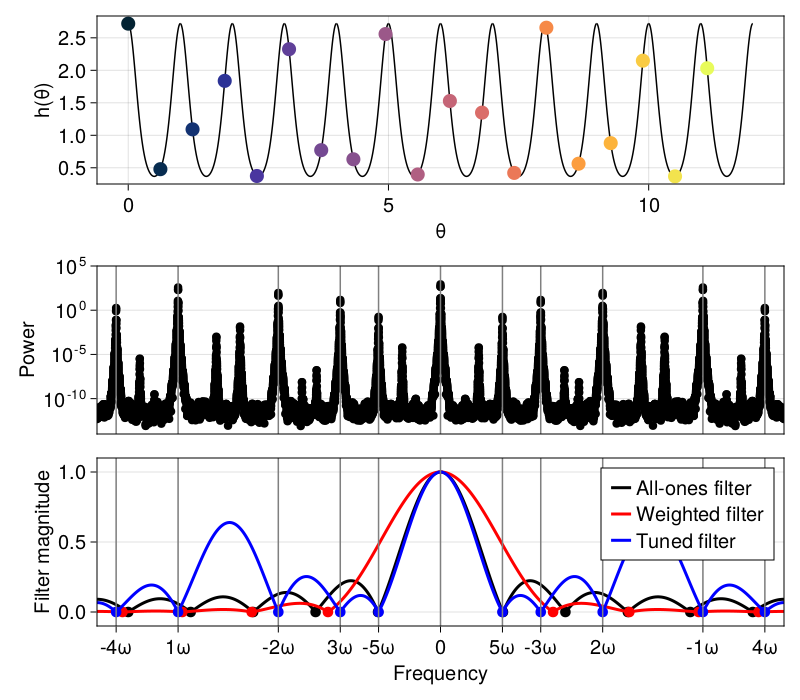}
    \caption{ (top) A test signal $\obsT(\theta) = e^{\cos 2\pi \theta}$, sampled at the points $\omega t$ where $\omega = (\sqrt 5 - 1)/2$. (middle) A discrete Fourier transform of the signal, showing peaks near the expected frequencies. (bottom) Three candidate $\Kback=11$ filters: the all-ones filter, a weighted Birkhoff filter, and an tuned filter to the first five roots of the frequency. We see that the all-ones filter is largest where there is a large amount of power if the signal, the weighted filter is small far from the zero frequency, and the tuned filter is zero exactly at the relevant frequencies. The absolute errors of the Birkhoff average for each filter are: (all-ones) $7.11 \times 10^{-2}$, (weighted Birkhoff) $7.38 \times 10^{-3}$, and (tuned) $2.72 \times 10^{-5}$.}
    \label{fig:background}
\end{figure*}

One way of understanding the weighted Birkhoff average is as a filter on the sequence $\abf_t = \obs(\sympmap^t(\xbf))$. When $\xbf$ is on an invariant torus, the sequence of observables has the form
\begin{equation}
\label{eq:a-sequence}
    \abf_t = \obsT(\omegabf t) = \sum_{n \in \Zbb} \obsT_{\nbf} \lambda_{\nbf}^t, \quad \lambda_{\nbf} = e^{2\pi i \nbf \cdot \omegabf}.
\end{equation}
So, the sequence $\abf_t$ is built out of equispaced samples of the trigonometric functions $\theta \mapsto e^{2\pi i \theta \nbf \cdot \omegabf}$. 
When we apply the weights $w_{k,\Kback}$ to the sequence and flip the order of summation, we find that
\begin{equation*}
    \WBA_{\Kback}[\obsT] = \sum_{n\in\Zbb} \obsT_{\nbf} q(\lambda_{\nbf}), \qquad q(\lambda) = \sum_{k = 0}^{\Kback-1} w_{k,\Kback} \lambda^k.
\end{equation*}
In this way, the weights act as a filter, where $w_{k,\Kback}$ will preferentially remove frequencies from the signal where the polynomial $q$ is small. An effective filter will return the mean $\obsT_0$, i.e.~$q(1) = 1$. The rest of the frequencies contribute to the error, so ideally $\abs{q(\lambda_{\nbf})} \ll 1$ for $\nbf \neq 0$. This is particularly important for the small $\abs{\nbf}$ modes that dominate the signal.

In Fig.~\ref{fig:background}, we see an example of this for an example observable $h(\theta) = e^{\cos (2\pi \theta)}$. 
We sample this signal at equispaced points $h(\omega t)$ where $\omega = (\sqrt 5-1)/2$ is the golden ratio. 
The top plot shows this signal and the samples. The middle plot shows a windowed discrete Fourier transform of the length $10000$ signal. 
The peaks of the power of the signal appear regularly at multiples of the frequency $\omega$, as is expected from the Fourier form of \eqref{eq:a-sequence}. 
Note that while the peaks occur at the locations $n\omega$, these peaks are not ordered as $n\omega$ wraps around the torus.
This is a consequence of the sampling being below the Nyquist sampling rate, a fact we do not have any practical control over.
In the bottom plot, we plot $\abs{q(e^{2\pi i \omega \theta})}$ for three potential filters with $\Kback = 11$:
\begin{itemize}
    \item The ``all-ones'' filter $w_{k,\Kback} = 1/\Kback$, used for the regular Birkhoff average.
    \item The weighted Birkhoff filter $w_{k,\Kback}$ sampled from the window function
\begin{equation*}
    w = e^{- ( t (1 - t))^{-1}}.
\end{equation*}
    \item A ``tuned'' filter that perfectly eliminates the first $\floor{\Kback/2}$ frequencies, found from the coefficients of the polynomial
    \begin{equation*}
        q_{\mathrm{tuned}} = \prod_{k=1}^{\floor{\Kback/2}}\frac{(z-e^{2\pi i \omega k})(z-e^{-2\pi i \omega k})}{(1-e^{2\pi i \omega k})(1-e^{-2\pi i \omega k})}.
    \end{equation*}
\end{itemize}
In each case, we can judge the absolute error of the filter's average by comparing it to the exact average (found via a weighted Birkhoff average with $\Kback=10^4$ to be $1.266066$). 
From Fig.~\ref{fig:background} (bottom), we see the all-ones filter polynomial has relatively large value at the peaks of the spectrum, resulting in the worst error of $7.11 \times 10^{-2}$. 
The weighted Birkhoff filter is much smaller at the spectral peaks with the most mass, resulting in a more accurate average with error $7.38 \times 10^{-3}$. 
The tuned filter does the best by two orders of magnitude, with an error of $2.72 \times 10^{-5}$. 

We note that while the tuned filter is small at the frequencies that dominate the signal, it is large in between. 
So, while the tuned filter worked well for this example, it would not work well if applied to a signal with a different value of $\omega$, as the polynomial roots are specifically related to multiples of the rotation number. 

\section{The Birkhoff Reduced Residual Extrapolation Method (Birkhoff RRE)}
\label{sec:method}
At the end of the previous section, we observed an important property: a filter that is tuned to specific frequencies in a signal can be significantly more effective than an arbitrary bump function. 
In this section, we will introduce a method to learn such a filter from a trajectory.


In Sec.~\ref{subsec:31}, we introduce a continuous problem for finding such an optimal filter on an ergodic component. Then, we discretize this problem with a Birkhoff average in Sec.~\ref{subsec:32}, which results in a variation of the reduced rank extrapolation (RRE) method \cite{Sidi2017}. Finally, in Sec.~\ref{subsec:33} we explain how we process the obtained filter to find the rotation number of invariant circles.

\subsection{The Least-Squares Problem}
\label{subsec:31}
We begin by defining a function for the action of a filter
\begin{equation*}
    \Fcal_{\Kback}[\obs](\xbf_0) = \sum_{k=0}^{\Kback-1} c_k (\obs \circ \sympmap^k)(\xbf_0),
\end{equation*}
and we call the associated filter polynomial $q_{\Kback}$. 
Note that this is equivalent to a regular Birkhoff average if $c_k = 1/\Kback$ and a weighted Birkhoff average if $c_k = w_{k,\Kback}$. 
Our goal is to find coefficients $c_k$ so that $\norm{\Fcal_{\Kback}[\obsT](\star) - \obsT_0}_{L^2}$ is small, where we are taking a norm over an ergodic region $X_0\subseteq X$ of the form
\begin{equation*}
    \norm{f}_{L^2}^2 = \int_{X_0} \abs{f(\xbf)}^2 \dif \mu.
\end{equation*}

However, there is a problem: we do not know $\obsT_0$ \textit{a priori}, so we cannot directly minimize $\norm{\Fcal_{\Kback}[\obsT](\star) - \obsT_0}_{L^2}$. 
So, we instead consider filtering the function 
\begin{equation*}
    \dobs(\xbf_0) = (\obs \circ \sympmap)(\xbf_0) - \obs(\xbf_0).
\end{equation*}
This is a convenient choice because $\dobs$ has zero mean and it is easy to calculate from a trajectory of $\obs$. Given $\dobs$, the new goal is to minimize $\norm{\Fcal \dobs(\cdot)}_{L^2}$, under the constraint that $\sum_k c_k = 1$. The constraint ensures the associated polynomial satisfies $q_{\Kback}(1) = 1$, a necessary condition to return the average. Additionally, on an invariant circle or island, $\dobs$ has the same type of Fourier series representation as $\obsT$. So, if a filter learns the rotation vector of $\dobs$, it will also be the rotation vector associated with $\obs$.

To discretize the norm, we use a weighted Birkhoff average \eqref{eq:WBA}: 
\begin{equation}
\label{eq:FKg}
    \norm{\Fcal_{\Kback} \dobs}^2_{L^2} \approx \WBA_T \left[\abs{(\Fcal_{\Kback} \dobs \circ \sympmap^t)(\star)}^2\right](\xbf_0).
\end{equation}
Assuming the conditions of Thm.~\ref{thm:Das2018} are met, we can bound the error of the above approximation by
\begin{equation*}
    \abs{\norm{\Fcal_{\Kback} \dobs}^2_{L^2(\Tbb)} - \WBA_T \left[\abs{(\Fcal_{\Kback} \dobs \circ R_\omega^t)(\star)}^2\right](\thetabf_0) } < C T^{-m}
\end{equation*}
for some $C > 0$ and integer $m$. 

In summary, the sum on the right hand side of \eqref{eq:FKg} can be seen as measure of how good a given filter $\filter$ is on an invariant set. 
Additionally, one can obtain this approximation using only a single trajectory of the dynamical system, rather than having \textit{a-priori} information. 
If we had used an unweighted Birkhoff average with no more assumptions, this returns the standard RRE algorithm. 
However, the weighted Birkhoff average acts as a convenient weighting for RRE by connecting it efficiently to a continuous problem. 

\subsection{Least Squares Solution}
\label{subsec:32}
Now that we have an energy to minimize, we discuss the numerical details. 
We begin by observing a structure of invariant circle and island signals: they come from pure Fourier series. 
That is, there are no growing or decaying modes, so the frequencies that we hope to learn via the filter all correspond to filter polynomial roots on the unit circle. 
Such roots do not change under time reversal (i.e. the conjugate pair obeys $(\lambda_{\nbf}, \bar{\lambda}_{\nbf}) \to (\bar{\lambda}_{\nbf}, \lambda_{\nbf})$ as $t \to -t$). This property corresponds to a linear constraint on the filter that 
\begin{equation}
\label{eq:time-reversal-invariance}
    c_{K+k} - c_{K-k} = 0,    
\end{equation}
where $0 \leq k < K$. Filters satisfying \eqref{eq:time-reversal-invariance} are known as \textit{palindromic}. We note that the converse is not true: time reversal symmetry does not imply roots on the unit circle. However, while not strictly necessary, we have found that this constraint dramatically improves the quality of the results. Throughout the rest of this paper, we will use $K$ to refer to a filter of length $2K+1$, whereas we used $\Kback$ and $T$ to represent filters of length $\Kback$ and $T$ respectively. We will find that $K$ is the number of unknowns of the final least-squares problem.

For the algorithm in this section, we assume that the user has access to an initial point $\xbf_0 \in X$, a symplectic map $\sympmap$, and an observable function $\obs$. 
The algorithm begins by sampling a trajectory of length $T+2K+1$ starting at $\xbf_0$ by repeated application of $\sympmap$. 
For many test maps, such as the standard map, this step is very quick. 
However, in applications where evaluating $\sympmap$ involves simulating a dynamical system up to a Poincar\'e section, this step could potentially dominate the cost of the algorithm.

From here, we compute the difference sequence $\ubf_t = \abf_{t+1} - \abf_{t}$ for $0 \leq t < T$. 
This step amounts to computing the difference function $\dobs \circ \sympmap^t$ from the previous section. 
Using this notation, the weighted Birkhoff average in \eqref{eq:FKg} can be expressed as the product
\begin{equation*}
    \WBA_T \left[\abs{(\Fcal_{\Kback} \dobs \circ \sympmap^t)(\star)}^2\right](\xbf_0) = \filter^T U^T W_T U \filter,
\end{equation*}
where $U \in \Rbb^{TD\times (2K+1)}$ is the block-Hankel matrix
\begin{equation*}
U = \begin{pmatrix}
        \ubf_{0}   & \ubf_{1} & \dots  & \ubf_{2K} \\
        \ubf_{1}   & \ubf_{2} & \dots  & \ubf_{2K+1} \\
        \vdots     & \vdots   & \ddots & \vdots\\
        \ubf_{T-1} & \ubf_{T} & \dots & \ubf_{T-1+2K}
    \end{pmatrix},
\end{equation*}
and $W_T \in \Rbb^{T D \times T D}$ is a diagonal matrix with the weighted Birkhoff weights
\begin{equation*}
   W_T = \begin{pmatrix}
        w_{0, T} I_D \\
        & \ddots \\
        & & w_{T-1, T} I_D
    \end{pmatrix},
\end{equation*}
where $I_D$ is the identity matrix in $\Rbb^{D\times D}$. 
We note that the matrix-vector product $U \filter$ can be interpreted as the filter being applied to sliding windows of $\ubf_t$, i.e.~$(U\filter)_t = \Fcal_K[\dobs\circ \sympmap^t](\xbf_0)$. 

For the full least squares problem, we require two more components. 
First, a filter must return the correct mean, which corresponds to the constraint that $\filter \cdot \bm{1} = q_K(1) = 1$. 
Second, we allow for a small amount of regularization to remove the possibility of low rank systems that arise from periodic orbits and very smooth circles. 
For the regularization, we choose weights so that the solution of the regularized problem is a weighted Birkhoff average. 
This is performed via the inverse of another weighted Birkhoff matrix $W_K \in \Rbb^{(2K+1) \times (2K+1)}$ with entries 
\begin{equation*}
    (W_K)_{kk} = \tilde{w}_{k,K} = \frac{w_{k, 2K+1} + w_{K-k, 2K+1}}{2}.
\end{equation*}
Note that $W_K$ has been symmetrized to respect the palindromic symmetry, and can equivalently be thought of as sampling the symmetric bump function $w(1/2 - x) + w(1/2 + x)$. 

In total, we define the following least-squares problem for finding $\filter$: 
\begin{gather}
\label{eq:lsqr1}
    R^2 = \min_{\filter \in \Rbb^{2K+1}} \norm{\begin{pmatrix}
        W_T^{1/2} U \\ \epsilon^{1/2} W_{K}^{-1/2}
    \end{pmatrix}\filter}^2,\\
\nonumber
    \text{s.t. } \bm{1} \cdot \filter = 1, \ \ c_{K+k} - c_{K-k} = 0 \text{ for } 1\leq k \leq K,
\end{gather}
where $\epsilon \geq 0$ is a regularization parameter. 
The above problem can be recognized as a weighted and time-reversal constrained version of RRE. 

We enforce the palindromic constraint by projecting by a matrix $P \in \Rbb^{K+1 \times 2K+1}$, where $\halffilter = P \filter$ with
\begin{align}
\label{eq:Pprojection}
    (P\filter)_k = \begin{cases}
        c_K, & k = 0, \\
        \frac{c_{K+k}+c_{K-k}}{\sqrt 2}, & 1 \leq K \leq K.
    \end{cases}
\end{align}
Substituting the matrix into \eqref{eq:lsqr1}, we find
\begin{gather}
\label{eq:LeastSquares}
R^2 = \min_{\halffilter \in \Rbb^{K+1}} \halffilter^T \tilde{A} \halffilter + \epsilon \halffilter^T \Whalf^{-1} \halffilter, \\
\nonumber
\text{s.t. } P\bm{1} \cdot \halffilter = 1, 
\end{gather}
where
\begin{equation*}
    \tilde{A} = P U^T W_T U P^T, \qquad \Whalf = P W_K P^T.
\end{equation*}
The effect of multiplying $U$ by $P^T$ is to ``fold $U$ in half'' and to sum the matching columns. 

Computationally, we can fully handle the constraints by including $\bm{1}\cdot \filter = 1$ in the projection. We do this by changing variables to a vector $\xibf\in\Rbb^{K}$ (see Sidi \cite{Sidi2017}, pp 41)
\begin{equation*}
    c_k = P_{\xi} \xibf + \ebf_K,
\end{equation*}
where $\ebf_K$ is the $K$th unit vector and
\begin{equation*}
     P_{\xi} = \begin{pmatrix}
        1 & -1 & & & & & & -1 & 1\\
         & \ddots & \ddots & & & & \iddots & \iddots  \\
         & & 1 & -1 & & -1 & 1 \\
         & & & 1 & -2 & 1
    \end{pmatrix}.
\end{equation*} 
Clearly, we have $P_{\xi} \bm{1} = 0$ and that $(P_{\xi})_{k,K+\ell} = (P_{\xi})_{k,K-\ell}$, so every constraint is satisfied. We can substitute this into our least-squares problem, giving
\begin{equation}
\label{eq:lsqr-practical}
    R^2 = \min_{\xibf \in \Rbb^{K}}  \norm{\begin{pmatrix}
        W_T^{1/2} U \\ \epsilon^{1/2} W_{K}^{-1/2}
    \end{pmatrix}P_\xi \xibf + \begin{pmatrix}
        W_T^{1/2} U \ebf_0 \\ \epsilon^{1/2} W_{K}^{-1/2} \ebf_0
    \end{pmatrix}}^2.
\end{equation}
We solve this system via a direct QR-based least squares solve, as this is reliably accurate and typically fast enough for the system sizes considered. 

We note that the system could alternatively be solved via an iterative least-squares solver such as LSQR \cite{paige1982}. The matrices $P_\xi$, $W_T^{1/2}$, and $W_{K}^{-1/2}$ can all be applied in $\Ocal(T)$ time. Additionally, the matrix $U$ can be applied $\Ocal(T\log T)$ time via fast Hankel multiplications using the fast Fourier transform. So, an iterative algorithm would have a worst-case run time $\Ocal(TK\log T)$ operations, rather than the $\Ocal(TK^2)$ time to perform the QR-based algorithm. Additionally, in the case that $K$ is chosen much larger than necessary, the iterative algorithm converges in many fewer steps. However, we have not yet found an appropriate preconditioner for this method: a necessary step for convergences to high tolerances.

In total, an algorithm for finding a filter is given below:
\begin{algorithm}[H]
\caption{Reduced Rank Extrapolation}
\label{alg:RRE}
\begin{algorithmic}[1]
\Require Initial point $\xbf_0$, symplectic map $\sympmap$, least squares dimensions $T$ and $K$, regularization $\epsilon$
\State Sample trajectory $\abf_t$ via repeated evaluations of $\sympmap$
\State Compute $\ubf_t$ and weights for $W_T$ and $W_K$
\State Solve \eqref{eq:lsqr-practical} for $\xibf$ via direct QR or iterative LSQR solver
\State Compute $\filter = P_\xi \xibf$ and $R$ from \eqref{eq:lsqr-practical}
\Ensure Best filter $\filter$ and residual $R$
\end{algorithmic}
\end{algorithm}

The above algorithm directly inherits convergence rates of the weighted Birkhoff average. 
\begin{theorem}
\label{thm:RRE-residual-convergence}
    Under Hypotheses \ref{setting}, the residual of \eqref{eq:lsqr1} converges independent of how $T$ depends on $K$ as
    \begin{equation}
    \label{eq:RRE-convergence-rate}
        \epsilon \leq R^2 \leq \epsilon + C_m (2K)^{-2m},
    \end{equation}
    where $m$ and $C_m$ are identical to Thm.~\ref{thm:Das2018}.
\end{theorem}
\begin{proof}
    Consider the solution $c_k = \tilde{w}_{k,K}$. By construction, the filter obeys the constraints. Using Thm.~\ref{thm:Das2018}, $\abs{(U\filter)_t} \leq C_m (2K+1)^{-m}$. So,
    \begin{equation*}
        \filter^T U^T W_T U \filter = \sum_{t=0}^{T-1} w_{t,T} (U\filter)_t^2 \leq C_m (2K)^{-2m}.
    \end{equation*}
    Similarly, the regularization term evaluates exactly as $\filter^T W_K^{-1} \filter = \epsilon$. Combining these computations, we have the upper bound.

    For the lower bound, we consider $\min_{\filter} \filter^T W_K^{-1} \filter \leq R$. Enforcing the $\bm{1} \cdot \filter = 1$ constraint with a Lagrange multiplier $\lambda$ and the palindromic constraint explicitly, the minimizer is found to solve the linear system
    \begin{equation*}
        \tilde{w}_{k,K}^{-1} c_k + \lambda = 0, \quad \sum_{k} c_k = 1.
    \end{equation*}
    This is solved via $c_k = \tilde{w}_{k,K}$ and $\lambda = -2$, the same solution as considered for the upper bound. Substituting $c_k$ back into the least-squares problem gives the lower bound.
\end{proof}
From \eqref{eq:RRE-convergence-rate}, we see that the extrapolation method converges at least as fast as the weighted Birkhoff average, particularly when $\epsilon = 0$. This means that we can also use the convergence of RRE to distinguish between chaotic and non-chaotic trajectories. 

In the case that $d=1$ and $\epsilon = 0$, the convergence rate can be improved using a continued fraction argument instead of relying on the weighted Birkhoff average:
\begin{theorem}
\label{thm:refpolyconvergence}
    Let $0 < \eta < 1$, $w : [0,1]\to \Rbb$ be a nonzero positive bounded function, and assume Hypothesis \ref{setting} (H3) with $d=1$ so that $\obsT \in C^M$. Then, for almost all $\omega\in\Tbb$, there exists an $L$ such that for $2K+1>L$ such that for some $C > 0$
    \begin{equation*}
        R^2 < C (2K+1)^{2\eta(-M+1)}.
    \end{equation*}
    Furthermore, if $\obsT \in C^\omega$ is a real analytic function on the torus, then for some $0 < r < 1$ independent of $\eta$ and $C > 0$, the error obeys the inequality
    \begin{equation*}
        R^2 < C r^{(2K+1)^{\eta}}.
    \end{equation*}
\end{theorem}
\begin{proof}
    See Appendix \ref{app:ideal-polynomial}.
\end{proof}
\begin{remark}
    We note that the proof of this theorem does not rely on any assumptions on $w$. That is, so long as $W_T$ are positive diagonal matrices with $\Tr(W_T) = 1$, the bound above will hold. Additionally, the nearly exponential rate obtained for analytic functions is faster than any known guarantee for the weighted Birkhoff average. This is due to the lack of reliance on bump functions, which can be $C^\infty$ but cannot be analytic.
\end{remark}

Because neither of the above theorems depend on $T$, there are not any strong lower bounds on the residual. 
Practically, it is required that $Td$ be greater than or equal to $K$ for a full rank linear system. 
This is an important requirement to ensure that chaotic trajectories do not converge. 
However, to our knowledge, there is currently no theorem here or in the weighted Birkhoff average literature proving that chaotic trajectories cannot be accelerated, despite strong numerical evidence. 
We will similarly give no proofs for the extrapolation method in chaos.

In order to classify trajectories efficiently, we also consider an adaptive algorithm. 
For this, we first define the scale-free residual $R_{G}$ as
\begin{equation}
\label{eq:RG}
    R_G = \frac{\sqrt{R^2-\epsilon}}{G}, \quad G^2 = \sum_{t = 0}^{T-1} w_{t,T} \abs{\ubf_t}^2 = \WBA[\abs{\dobs}^2](\xbf_0).
\end{equation}
Because $G$ converges to a nonzero number with $T$ (assuming $\xbf_0$ is not a fixed point), $R_G$ converges at the same rate as $R$. 
For the adaptive algorithm, we increase $K$ from an initial value $K_{\min}$ to a maximum value $K_{\max}$ by an increment of $\Delta K$.
We assume that $T$ scales with $K$ as $T = \ceil{\gamma K/D}$ for some constant $\gamma \geq 1$. 
At each increment, we check whether $R_G < \delta$, and finish the algorithm when this conditions is met.
The fact that $R_G$ is scale-free allows for the algorithm to be applied to multiple orbits for the same map with an accuracy to match the largest scale of the system.
In total, the adaptive algorithm is below:
\begin{algorithm}[H]
\caption{Adaptive Reduced Rank Extrapolation}
\label{alg:adaptiveRRE}
\begin{algorithmic}[1]
\Require Initial point $\xbf_0$, symplectic map $\sympmap$, regularization $\epsilon$, $T$ proportionality constant $\gamma$, convergence cutoff $\delta$, initial $K$ value $K_{\mathrm{init}}$, maximum $K$ value $K_{\max}$, $K$ increment $\Delta K$
\State $K \gets K_{\mathrm{init}}$, $T \gets \ceil{\gamma K/D}$
\While {$K \leq K_{\max}$ and $R_G > \delta$}
    \State Get $\filter$, $R$ via Algorithm \ref{alg:RRE}, reusing trajectory information
    \State Get $R_G$ via \eqref{eq:RG}
    \State $K \gets K + \Delta K$, $T \gets \ceil{\gamma K/D}$
\EndWhile
\Ensure Best filter $\filter$ and residual $R$
\end{algorithmic}
\end{algorithm}
With the use of $QR$ factorization updates, the adaptive algorithm maintains an $\Ocal(T K^2)$ runtime.
However, even without updates, choosing large enough steps $\Delta K$ mostly avoids the cost of adaptivity even without updates for reasonably small systems. 

\subsection{Processing the Learned Filter}
\label{subsec:33}
Once we have found a filter $\filter$ for the sequence of observations for a trajectory, we must post-process the results. 
There are two steps to this process. 
The first is straightforward: we distinguish chaotic trajectories from non-chaotic trajectories via a tolerance for the residual. 
The same tolerance as in Algorithm \ref{alg:adaptiveRRE} can be used, although often it is convenient to use separate tolerances for adaptivity convergence and chaos classification. 
If the residual is above the tolerance, we deem the trajectory ``chaotic,'' and otherwise it is an invariant circle or an island. 
We note that choosing the tolerance for a specific problem depends on the desired accuracy, the tolerable amount of computation, and the accuracy of the $\sympmap$ evaluation. 
However, any error in the computation of $\filter$ will lead to more difficulty in the following steps, so lower tolerances are always more reliable for adaptive RRE than classification.

In the case that we deem the trajectory non-chaotic, the second step is to process the roots of the filter polynomial.
Empirically, we have found the following convergence result to hold for the polynomial roots:
\begin{conjecture}
\label{con:least-squares}
    Under Hypotheses \ref{setting}, let $\filter_K$ be a sequence of solutions to \eqref{eq:LeastSquares} with $T = \ceil{\gamma K}$ for $\gamma \geq 1$ and let $q_K(z)$ be the filter polynomial with coefficients $\filter_K$. Then, for all $n$ such that $\abs{\dobsT_n} \neq 0$, there exists an $M_* > 0$ and a sequence $z_{K,n}$ such that $q_K(z_{K,n}) = 0$ and
    \begin{equation}
    \label{eq:lambdaconvergence}
        \abs{\lambda_n - z_{K,n}} < C_n K^{-M_*},
    \end{equation}
    where $M_*$ depends on $M$ and potentially $\nu$. 
    In particular, if $M=\infty$, $M_*$ is unbounded.
\end{conjecture}
The above conjecture essentially states that the roots of the filter necessarily approach the points $\lambda_n$ associated with multiples of the rotation number and the island period. 
The difficulty with in proving the conjecture is that the weighted Birkhoff average performs strongly for minimizing the RRE residual.
So, to prove the roots converge, one must show that filter polynomials with roots that converge to $\lambda_n$ out-perform weighted Birkhoff averages in minimizing the residual.

Numerically, roots of the filter polynomial are found by solving the equation
\begin{equation*}
    q_K(\lambda) = \sum_{k = 0}^{2K} c_k \lambda^k = 0,
\end{equation*}
for every root $\tilde \lambda_j$. 
This is performed by reducing the polynomial to another Chebyshev polynomial with half the dimension, then solving for the eigenvalues of the associated colleague matrix \cite{good1961} (see App.~\ref{app:palindromic-roots}). 
The dimension reduction accelerates the root-finding by a significant constant factor.
Additionally, due to the majority of eigenvalues being on the unit circle, this eigenvalue problem is well conditioned. 
The only exception to this is when $\lambda = -1$ is a root, as happens with $p=2$ islands. 
In this case, the palindromic coefficients force $-1$ to be a double root, causing a square root of the error of this frequency. 

Throughout the rest of this section, we assume that $d=1$ and therefore the trajectory is either an invariant circle or a $d=1$ island.
Using the roots $\tilde \lambda_n$, we will find the rotation number by finding which roots represent the majority of the signal.
This step is heuristic: we assume that the low-frequency oscillations will dominate the Fourier spectrum, so we can use this information to prune the important frequencies from the less important ones. 
To do this, we first restrict to values of $\tilde \lambda_n$ that are very close to the unit circle by some tolerance on the order of the square root of machine epsilon (due to the aforementioned sensitivity of $-1$ as a root). 
The remaining eigenvalues are sorted by solving the weighted least-squares system for each eigenvalue's prominence in the signal
\begin{equation*}
    \min_{V} \norm{W_{2K+T+1}^{1/2} \left( \tilde \Phi V  - A\right)}^2,
\end{equation*}
where $\tilde \Phi \in \Rbb^{2K+T+1 \times 2K+1}$ with $\tilde \Phi_{mn} = \tilde \lambda_{n}^{m}$ is the matrix of eigenmodes associated with the frequencies $\tilde \lambda_n$, $V \in \Rbb^{2K+1 \times D}$ holds the prominence of each mode in the signal, $W_{2K+T+1}\in \Rbb^{2K+T+1\times2K+T+1}$ is the weighted Birkhoff diagonal matrix with $(W_{2K+T+1})_{mm} = w_{m,2K+T+1}$, and 
\begin{equation*}
    A = \begin{pmatrix}
        \abf_0^T \\ \vdots \\ \abf_{T+2K}^T
    \end{pmatrix}.
\end{equation*}
The eigenvalues are then sorted by the value of norm of the rows of $V$.

Using the top few eigenvalues (arbitrarily chosen to be $10$), we next determine if the sequence corresponds to an island. 
We do this by noticing a key difference between a $p\geq 2$ island and a Diophantine invariant circle: the island spectrum \eqref{eq:a-island} has rational frequencies whereas the invariant circle \eqref{eq:a-sequence} does not. 
So, if the sequence is an island, we expect a prominent rational frequency $\tilde \omega_n = (2\pi i)^{-1} \arg \tilde \lambda_n$ in the sorted eigenvalues. 
For this, we use a Farey-sequence method from Sander and Meiss\cite{sander2020} (Algorithm 2) to determine whether a given root is rational to a specified tolerance $\epsilon_{\mathrm{rat}}$ and bounded denominator $p \leq p_{\max}$. 
If one of the roots is rational by this algorithm, we consider it to correspond to an island chain with the largest found period $p$. 
From here, we can rerun Algorithm \ref{alg:RRE} on the ``stacked'' signal with dimension $\hat{D} = p D$
\begin{equation*}
    \hat{\abf}_t = \begin{pmatrix}
        \abf_{pt} \\ \abf_{pt+1} \\ \vdots \\ \abf_{pt+p-1}
    \end{pmatrix},
\end{equation*}
with a filter length $\hat{K} = \floor{K/D}$.
Using this signal is equivalent to working simultaneously with the $p$ invariant circles of $\sympmap^p$ with the same rotation number. 
The resulting filter will be that of an invariant circle and the above steps can be repeated, skipping the check for rational roots.

Now, consider the case that no rational roots of the filter are found. 
In order to determine the rotation number, we simply choose the frequency that appears first in the sorted eigenvalues, i.e. $\omega = (2\pi i)^{-1} \arg(\tilde \lambda_0)$ where $\tilde \lambda_0$ is the first eigenvalue after sorting. 
We note that this is not always true: invariant circles can be found where the largest Fourier coefficient does not correspond to the rotation number. 
However, we have not found this to be a practical issue for any of our examples.

Once the rotational frequency $\omega$ is determined, we can obtain the coefficients of the invariant circle associated with the signal by solving another linear system. This time, because typically $\omega$ is resolved to very high accuracy, we solve the linear system
\begin{equation}
\label{eq:circle-projection}
     \min_V \norm{W_{2K+T+1}^{1/2} \left( \Phi V - A\right)}^2,
\end{equation}
where $\Phi \in \Rbb^{2K+T+1 \times L}$ with $\Phi_{mn} = \lambda_n^m$, $\lambda_n = e^{2\pi i \omega (n-L)}$ are the new signal frequencies, 
and $L$ is the number of coefficients wanted, chosen to control the condition number of the least-squares problem (see Sec.~\ref{app:choosing-L}).
The approximated invariant circles are then given as
\begin{equation*}
    \zbf^{(j)}(\theta) = \sum_{\ell = -L}^L V_{\ell,D(j-1):Dj-1}e^{2\pi i \ell \theta},
\end{equation*}
where we are using ``Matlab notation'' for slicing the rows of $V$ into each island component. 

Once an invariant circle or island chain is found, the fit can be validated via a residual from the parameterization. 
In particular, we can evaluate the residual
\begin{align}
\nonumber
    R_p^2 &= \frac{1}{pJ} \sum_{j=0}^{J-1} \left(\vphantom{\sum_{j=1}^{p-1}}\abs{\zbf^{(1)}(j h + \omega) - \sympmap(\zbf^{(p)}(j h))}^2\right. \\
    &\quad \qquad \qquad +  \left.\sum_{j = 1}^{p-1} \abs{\zbf_{(j+1)}(hj) - \sympmap(\zbf^{(j)}(hj))}^2 \right), \\
\label{eq:validation-error}
    &\approx \frac{1}{p} \left(\vphantom{\sum_{j=1}^{p-1}}\int_{\Tbb}\abs{\zbf^{(1)}(\theta + \omega) - \sympmap(\zbf^{(p)}(\theta))}^2 \dif \theta\right.\\
    & \quad \qquad \qquad + \left.\sum_{j = 1}^{p-1} \int_{\Tbb}\abs{\zbf_{(j+1)}(\theta) - \sympmap(\zbf^{(j)}(\theta))}^2 \dif \theta \right),
\end{align}
where $h = 1/J$ for $J \in \Nbb$. The residual $R_p$ is essentially an $L^2$ measurement of the conjugacy \eqref{eq:island-conjugacy}. If this residual is small, then the conjugacy is likely correct. Furthermore, a small $R_p$ likely puts the island within the basin of convergence for the parameterization method \cite{haro2016}, which could be used to refine the estimate of the island. Thus, Birkhoff RRE could be seen as a method of finding an initial guess for higher-accuracy methods.
\section{Examples}
\label{sec:examples}

\begin{figure*}
    \centering
    \includegraphics[width=0.7\textwidth]{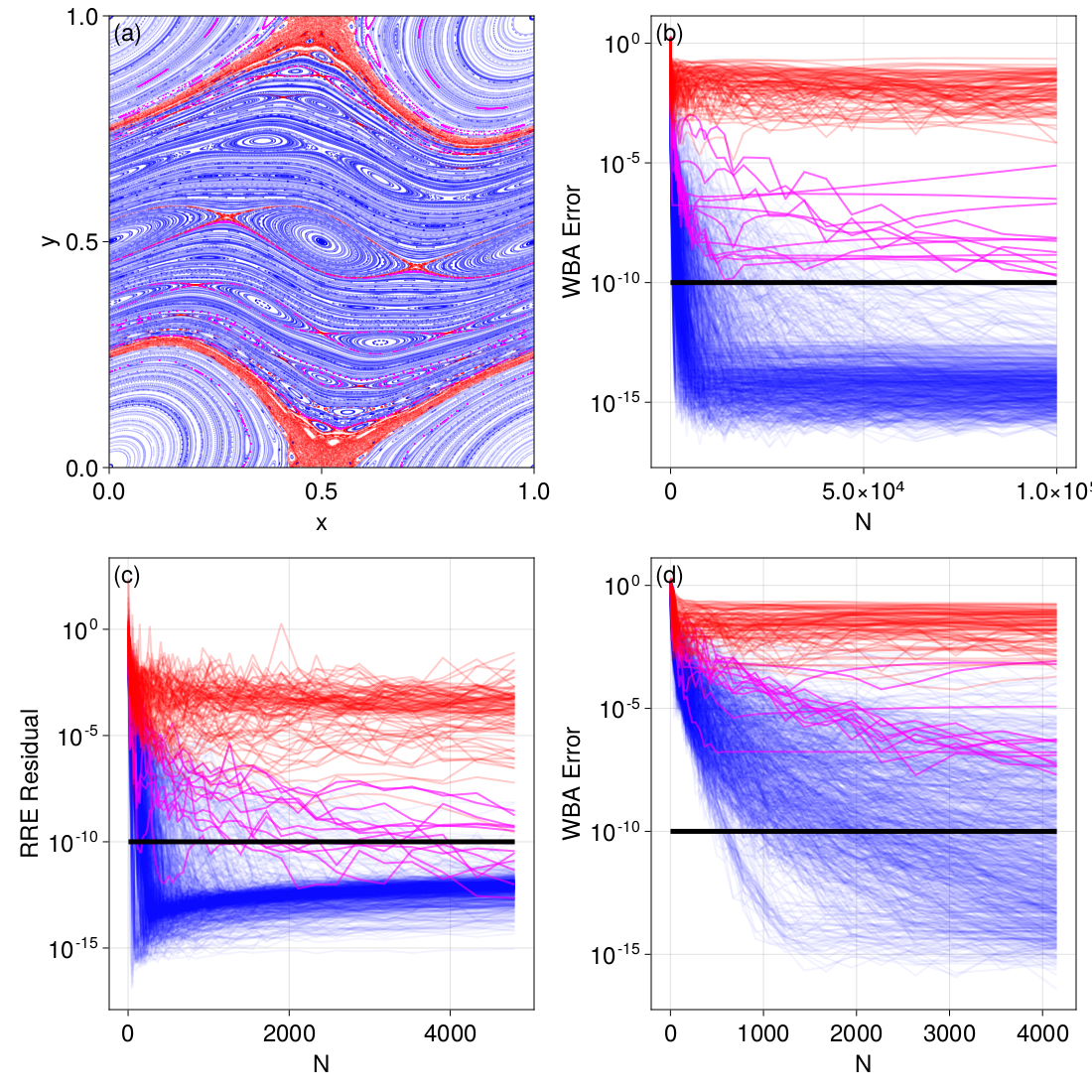}
    \caption{(a) A Poincar\'e plot of the standard map with $k=0.7$, colored by trajectory classification of integrable, chaotic, or indeterminate. (b) Convergence of the weighted Birkhoff doubling error with respect to trajectory length. Trajectories are colored by the ending point on this plot, with $R > 10^{-5}$ chaotic, $R < 10^{-11}$ integrable, and the rest indeterminate. (c) Convergence of RRE residual with respect to trajectory length. (d) The same as data as (b), on the same domain as (c). From (c) and (d), we see that the RRE residual appears to converge much more rapidly than the weighted Birkhoff average. }
    \label{fig:residual-convergence}
\end{figure*}

In this section, we give two methods of the Birkhoff RRE method applied to symplectic maps. 
In Sec.~\ref{subsec:standard-map}, we apply Birkhoff RRE to the Chirikov standard map in order to investigate its convergence properties.
We show that the convergence of the RRE residual is significantly more efficient in the number of samples than the weighted Birkhoff average (WBA).
As a consequence, RRE can classify trajectories with many fewer symplectic map iterations.
We also show that the roots of the filter polynomial converge like Conjecture~\ref{con:least-squares}.

In the Sec.~\ref{subsec:stellarator}, we explore the performance of Birkhoff RRE on magnetic field-line dynamics for a stellarator, a type of plasma confinement device.
We show how the method classifies chaotic and integrable trajectories, as well as finding Fourier representations of circles and islands.
The rotation number and residuals are also reported.
The symplectic map is obtained by numerically integrating magnetic field lines from an interpolated field, so this example additionally shows the performance on a map with symplecticity breaking error.

Code for performing both Birkhoff RRE and the weighted Birkhoff average herein can be found in the \texttt{SymplecticMapTools.jl}\footnote{\texttt{https://github.com/maxeruth/SymplecticMapTools.jl}} Julia package.

\subsection{Standard Map Convergence}
\label{subsec:standard-map}
For the first experiment, we are interested in comparing the convergence and classification of Birkhoff RRE vs WBA. 
For RRE, we measure convergence via the square root of the least-squares residual $R$ defined in \eqref{eq:LeastSquares}. 
For WBA, we use the method of Sander and Meiss\cite{sander2020}, and compare the values of two consecutive averages. 
That is, given a symplectic map $\sympmap$ and an observable $\obs$, the WBA residual is
\begin{align*}
    R_{\text{WBA}} &= \norm{\WBA_T[\obs](\xbf_0) - \WBA_T[\obs](\sympmap^T (\xbf_0)) },\\
    &= \norm{\sum_{t=0}^{T-1}w_{t,T}\abf_t - \sum_{t=0}^{T-1}w_{t,T} \abf_{t+T} }.
\end{align*}
Both residuals have the same theoretical rate guarantees for $C^\infty$ functions, so it is necessary to numerically check that RRE converges significantly faster.

For a fair comparison, we hold the total number of map evaluations $N$ constant between the two methods. 
In particular, we compare $N = (2+\gamma)K+1$ for RRE against $N = 2T$ for WBA, where we choose the ``rectangularity'' constant to be $\gamma = 2$. 
The example dynamical system we compute with is the Chirikov standard map \eqref{eq:standard-map} with $\ksm=0.7$ (see also Fig.~\ref{fig:standard-map}). In order to handle the map domain $\Tbb \times \Rbb$, we use the smooth observable $\obs : \Tbb \times \Rbb \to \Rbb^2$
\begin{equation*}
    \obs(x,y) = (y+0.5)\begin{pmatrix}
        \cos(2\pi x) \\ \sin(2\pi x)
    \end{pmatrix}.
\end{equation*}

\begin{figure}
    \centering
    \includegraphics[width=0.4\textwidth]{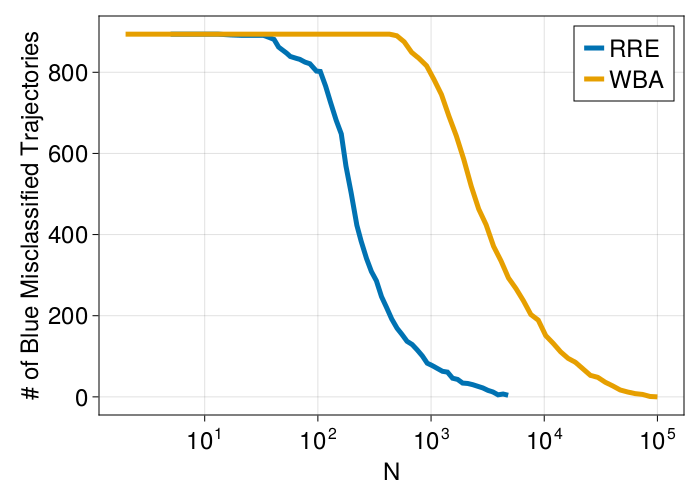}
    \caption{Number of integrable trajectories of Fig.~\ref{fig:residual-convergence} misclassified, with the tolerance of both WBA and RRE set to $10^{-11}$. We see that RRE converges to a low misclassification rate much more efficiently than WBA in the number of map iterations.}
    \label{fig:misclassification}
\end{figure}

\begin{figure*}
    \centering
    \includegraphics[width=0.8\textwidth]{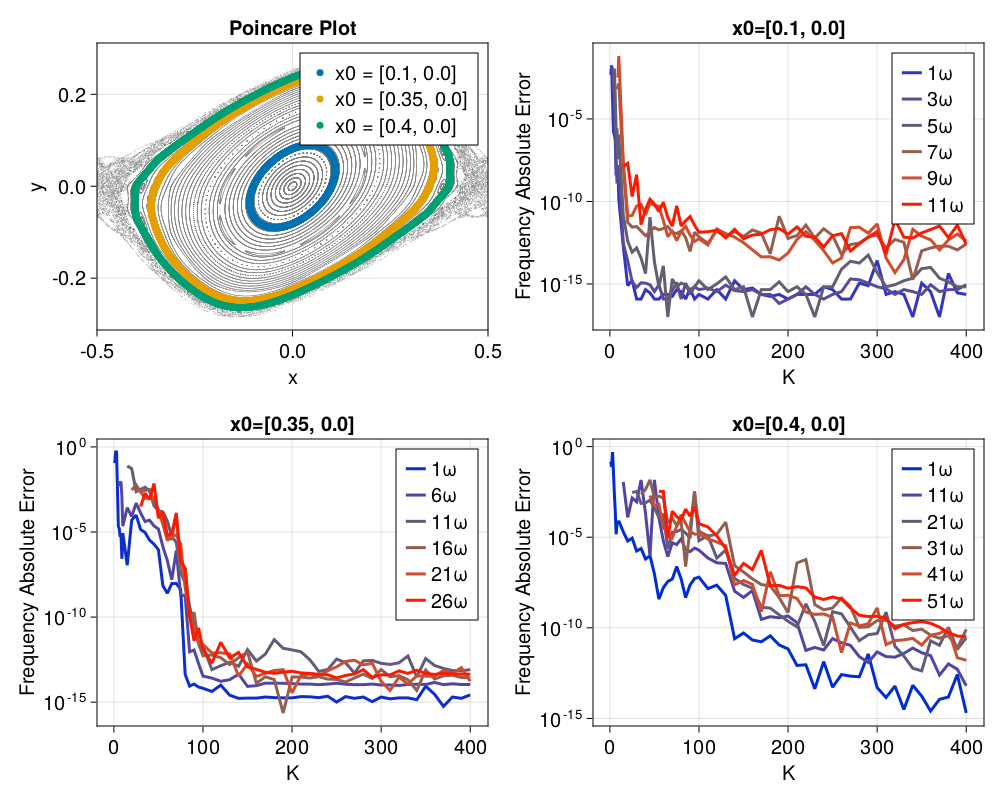}
    \caption{(a) A Poincar\'e plot of the standard map, with three invariant circle with initial $x$ points $0.1$, $0.35$, and $0.4$. (b-d) Convergence plots for the absolute error of the learned frequencies from the polynomial filter vs filter length $K$. The absolute error is computed by comparing to the frequency learned from a $K=1500$ simulation. In (b), we see very rapid convergence for the smooth central circle. In (c), we see slower, but still rapid convergence of the learned frequencies. In (d), for an invariant circle near the edge of chaos, we see significantly slower convergence. The errors are approximately straight in this plot, indicating nealry exponential convergence of the roots.}
    \label{fig:frequency-convergence}
\end{figure*}

A Poincar\'e plot of the standard map is found in Fig.~\ref{fig:residual-convergence} (a), where $1000$ trajectories are plotted. 
The trajectories are classified by a long-time $N=10^5$ weighted Birkhoff average, with the three categories of nested invariant circles and islands (blue, $R_{\text{WBA}} < 10^{-11}$), chaos (red, $R_{\text{WBA}} > 10^{-5}$), and indeterminate (purple, $10^{-5} \leq R_{\text{WBA}} \leq 10^{-11}$. 
The convergence of the WBA residual used to classify these trajectories as a function of $N$ is found in Fig.~\ref{fig:residual-convergence} (b).
Visually, the `indeterminate' trajectories are typically either at the transition of trajectory types (e.g.~circles to islands or chaos) or have a near small-denominator rotation number.
We note that difficult classification on transitional domains is the typical behavior of any numerical method depending on a continuous quantity.
Additionally, misclassifications are likely to happen for any trajectory classification algorithm, as orbits on one side of a classification boundary can shadow orbits on the other side for arbitrarily long periods of time. 

In Fig.~\ref{fig:residual-convergence} (c), we plot the convergence of the RRE residual for $N \leq 2801$ and $\epsilon = 0$.
For comparison, we plot the WBA residual on the same domain in Fig.~\ref{fig:residual-convergence} (d). 
We see that RRE converges significantly quicker than WBA, with most residuals reaching the machine precision limit before $N=1000$. 

In Fig.~\ref{fig:misclassification}, we compare the number of misclassified integrable trajectories for RRE and WBA.
To obtain this number, we subtract the number of blue trajectories below the $10^{-11}$ tolerance (shown as a black horizontal line in Figs.~\ref{fig:residual-convergence} (b)-(d)) from the total number of blue trajectories for each value of $N$.
We see that RRE classification happens an order of magnitude faster than that of WBA.

Next, we consider the convergence of the roots of the filter polynomials. 
For this, we choose three invariant circles of the standard map of varying smoothness, plotted in Fig.~\ref{fig:frequency-convergence} (a).
The inner circle is near the core of the nested circles, and is well approximated by a small number of Fourier modes. 
The yellow circle is more complex, and the green circle is a case on the edge of chaos. 

For each circle, we perform RRE to high accuracy ($K=1500$, $\gamma = 3$) and obtain the rotation number $\omega$ by the process described in Sec.~\ref{subsec:33}. Then, for increasing $K$, we let $q_K$ be the filter polynomial for $\gamma = 3$ and let $\zbf_{K,n}$ be the associated roots. From Conjecture \ref{con:least-squares}, we expect the roots to converge faster than $\Ocal(K^{-M_*})$ for all $M_*$. To measure this, in Fig.~\ref{fig:frequency-convergence} (b-d) we plot the error
\begin{equation*}
    E_{K,m} = \min_{n} \abs{\zbf_{K,n} - e^{2\pi i m \omega}}
\end{equation*}
for a variety of values of $m$ for each circle. In all cases, we find that the convergence for lower multiples of the rotation number is faster than that for higher multiples. This likely due to the higher prominence of these Fourier modes in the signal. For both the small and medium circles, the rotation number converges to the correct value for values of $K$ less than $100$, corresponding to less than $500$ iterations of the map. For the outer circle, the roots converge significantly more slowly, reaching machine precision at around $K=400$. Additionally, the lines are approximately straight, indicating a nearly exponential convergence of $\Ocal(e^{-\alpha K})$ for some $\alpha > 0$. As the circles become more difficult to approximate, higher multiples of the rotation number are found to converge due to a higher number of modes with nontrivial Fourier coefficients, with over $50$ multiples converging to high accuracy for the outer circle.

\subsection{A Stellarator Example}
\label{subsec:stellarator}


\begin{figure*}
    \centering
    \includegraphics[width=0.9\textwidth]{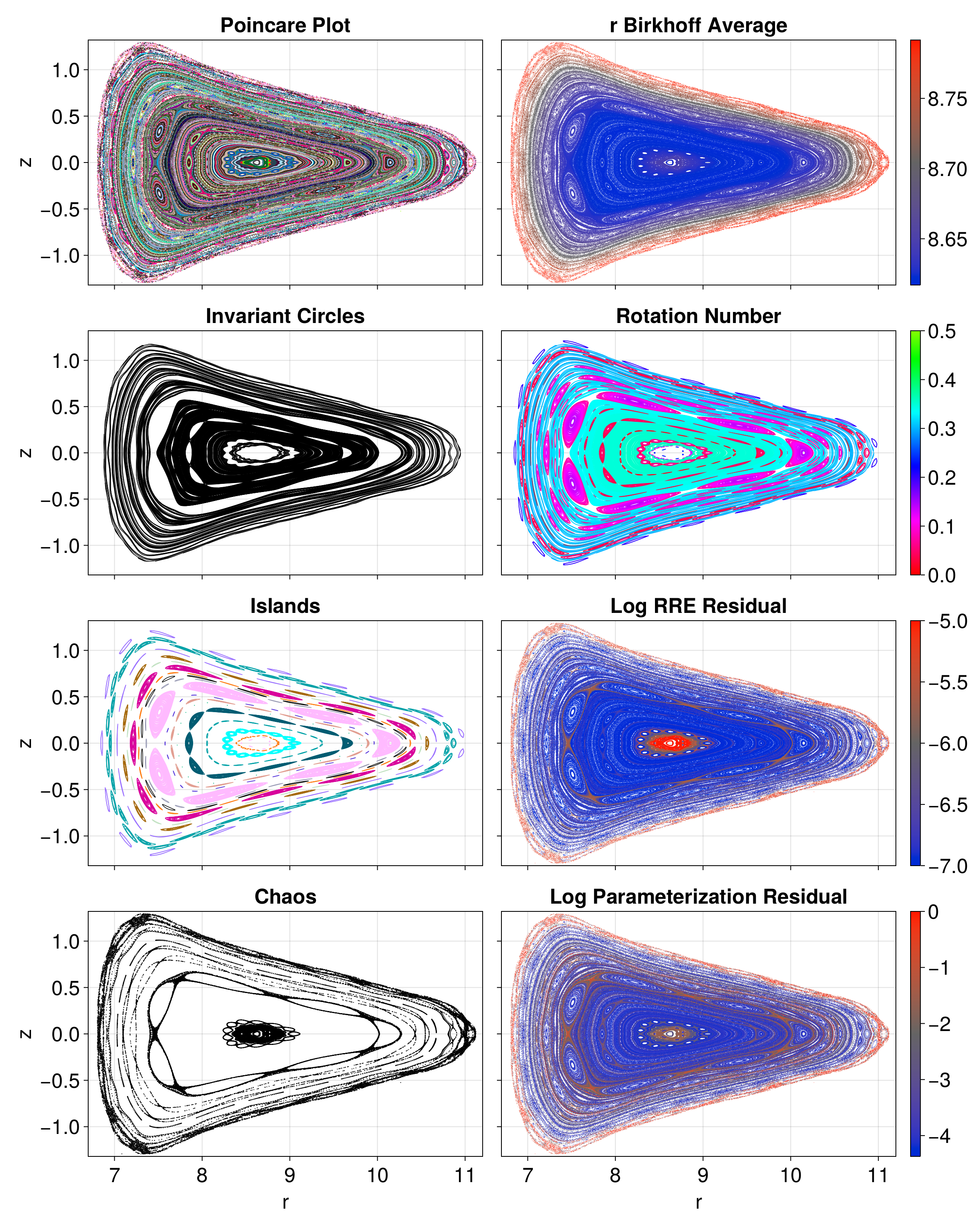}
    \caption{Eight plots of an optimized stellarator. On the left, we show a Poincar\'e plot and the trajectories separated into invariant circles, islands, and chaos. On the right, we show the Birkhoff average of the $r$ coordinate, the rotation number, the log RRE residual, and the log parameterization residual, all computed using the Birkhoff RRE procedure.}
    \label{fig:stell-all-plots}
\end{figure*}

In this section, we consider our method applied to stellarators, a type of toroidal plasma confinement device. 
Stellarators use large magnetic fields to confine the charged particles that comprise a plasma within the device.
In particular, high energy particles approximately follow magnetic field lines, described by the dynamical system
\begin{equation*}
    \dot{\xbf} = \Bbf(\xbf),
\end{equation*}
where $\Bbf$ is the magnetic field and $\xbf = r \ebf_r(\phi) + z \ebf_z$ is the cylindrical position vector with radial position $r$, azimuthal angle $\phi$, and $z$ is the vertical position. 
Magnetic field line dynamics is a 1.5D Hamiltonian system, and a symplectic map is obtained by numerically integrating the magnetic ODE over a field period in $\phi$. 
We note that the stellarator map has another difference from the standard map: the domain $(r, z) \in \mathbb{R}^2$ does not have the rotation number ``built in.'' 
In the case of the standard map, an average of the $y$ coordinate gives the rotation number, whereas the magnetic axis of the stellarator is not necessarily known \textit{a priori}.

To evolve the magnetic field, we use the Julia \texttt{Tsit5} integrator on a configuration with many islands.
The magnetic field is interpolated from the output of a plasma equilibrium solver using cubic splines. 
A Poincar\'e plot showing the $1000$ trajectories is found in Fig.~\ref{fig:stell-all-plots}(a).
Due to both the integrator and the field interpolation, there are small non-symplectic errors in the map.
These errors affect the clarity of classifying trajectories, as the Birkhoff RRE residuals do not reach the same low levels as the previous example.
The errors also make identifying rational numbers and rotation numbers more difficult, and the associated tolerances must also be adjusted. 
This potentially leads to a increased chance of misclassification, so we perform the following simulation to ensure Birkhoff RRE is robust to errors.

We perform the adaptive Birkhoff RRE classification algorithm \ref{alg:adaptiveRRE} on $1000$ trajectories using $\obs : (r,z) \mapsto (r,z)$, $K_{\min} = 50$, $K_{\max} = 400$, $\Delta K = 50$, $\gamma = 3$, $\epsilon = 0$, and $\delta = 10^{-7}$. 
With these parameters, a maximum of $(2+\gamma)K_{\max}+1 = 2001$ iterations of the map are used per invariant circle, and a minimum of $(2+\gamma)K_{\min}+1 = 251$ iterations are used. 
We note that the cutoff $\delta$ is chosen low enough to distinguish chaos, but high enough that it is resilient to the noise of evaluating the return map (compare to the value of $10^{-11}$ for the standard map). 
After the adaptive RRE algorithm is run, we attempt to fit an invariant circle or island chain to every trajectory, regardless of whether the algorithm converged before the tolerance. 
Then, we evaluate the validation error \ref{eq:validation-error} for each invariant circle.

We plot the invariant circles and islands parameterizations in Fig.~\ref{fig:stell-all-plots}(b-c), with the islands colored by the period $p$. 
We only plot those with a validation error below the threshold $R_p < 10^{-2}$. 
The associated rotation numbers to the circles and islands are plotted in Fig.~\ref{fig:stell-all-plots}(f). 
In Fig.~\ref{fig:stell-all-plots}(d), we plot the chaotic trajectories, defined as those with the RRE residual satisfying $R > 5 \times 10^{-7}$. 
We note that this value is different than the tolerance used for classification. 
This is because the standard needed to classify as integrable is separate from the accuracy needed to obtain an effective parameterization. 

The circles, islands, and chaos all match expectations from the Poincar\'e plot.
The classification algorithm makes it clear that much of the volume is taken by nested invariant circles, with many islands of various periods intermixed. 
The classification successfully identifies high-period and low-radius islands that may have been missed without computer processing. 
Additionally, chaotic trajectories are found at the boundaries of invariant circles and islands, matching the intuition provided by the standard map example. 
The rotation number plot gives a visual representation of the shear.
This allows us to identify the island chains with specific resonant low-denominator frequencies, e.g.~with the largest islands occurring at $\omega = 2/7$, $3/11$, and $3/10$. 

In Fig.~\ref{fig:stell-all-plots}(f-g), we plot the Birkhoff RRE residual $R$ and validation error $R_p$ on a logarithmic scale. 
We find the two plots have similar features: both residuals tend to be larger in regions of chaos (the inner core and outer ring) and smaller where there are integrable trajectories. However, we find that the parameterization residual tends to give a sharper indicator of correctness, with chaotic trajectories typically appearing as more red.

\section{Conclusions}
\label{sec:conclusions}
We have shown how Birkhoff RRE, a version of the reduced rank extrapolation algorithm, can be used to accelerate ergodic averages on invariant tori.
We find numerically find that this acceleration is significantly more stronger than the weighted Birkhoff average for the same number of map evaluations, at the cost of a more complex algorithm.
The acceleration can be used to efficiently classify trajectories in symplectic maps as integrable or chaotic.
Beyond this, using the resultant filter from, we can compute the number of islands and the rotation number of integrable trajectories to high accuracy.

Due to Birkhoff RRE post-processing step's reliance on Conjecture \ref{con:least-squares}, a proof would be a natural future direction from this work.
The proof of this convergence would contrast with standard convergence proofs for RRE, which all rely all but finitely many of the frequencies $\lambda_j$ having modulus less than one. 
This does not hold for Birkhoff RRE (we have $\abs{\lambda_j} = 1$ for all $j$), so the proof would likely instead rely on a combination of the Diophantine property and smoothness of $\obsT$, just as both the KAM theorem and the weighted Birkhoff average do.

Beyond this, it would be interesting to apply Birkhoff RRE to higher-dimensional tori.
We have proven that the Birkhoff RRE residual converges in higher dimensions, so it should be able to effectively classify trajectories.
Additionally, it is likely that the frequencies converge with similar rates, allowing for high-dimensional Diophantine vectors to be identified without any winding argument.
This would be a convenient option for processing high-dimensional symplectic geometry without the need for complicated initial guesses or continuation.

\appendix

\section{Finding the Roots of a Palindromic Polynomial}
\label{app:palindromic-roots}
Here, we detail how we solve for the roots of palindromic polynomials. This method is no more or less accurate than solving for the eigenvalues of the standard Frobenius companion matrix (see the book by Trefethen\cite{trefethen2013} for an introduction). However, this process halves the dimension of the eigenvalue problem, resulting a significantly faster method than the standard companion matrix approach. 

The algorithm results from the observation that palindromic polynomials of the form 
\begin{equation*}
    P(z) = \sum_{k = 0}^{2K} c_k z^{k},
\end{equation*}
can be reduced to degree $K$ polynomials in $(z+z^{-1})/2$. In particular, we have that
\begin{equation*}
    z^{-K} P(z) = Q\left(\frac{z+z^{-1}}{2}\right), \qquad Q(x) = \sum_{k=0}^K b_k T_k(x),
\end{equation*}
where $T_k$ are the Chebyshev polynomials. If all of the roots of $P$ are on the unit circle, the change of variables implies that all the roots of the $Q$ are on the interval $[-1,1]$, so we expect Chebyshev polynomials to be well conditioned for this problem. Typically results of Birkhoff RRE have almost all of the roots on the unit circle for non-chaotic orbits, so they satisfy this assumption. Additionally, using the Chebyshev three-term recurrence
\begin{align*}
    T_0(x) &= 1, \\
    T_1(x) &= x, \\
    T_{n+1}(x) &= 2x T_n - T_{n-1},
\end{align*}
one sees that the polynomials in $z$ coordinates have a convenient form of
\begin{equation*}
    T_k\left(\frac{z+z^{-1}}{2} \right) = \frac{z^k + z^{-k}}{2}.
\end{equation*}
This results in the convenient relations
\begin{equation*}
    b_k = \begin{cases}
        c_K, & k=0, \\
        2c_{K+k}, & 0<k\leq K.
    \end{cases}
\end{equation*}

To find the roots of $p$, we use the Chebyshev colleague matrix
\begin{equation*}
    C = \begin{pmatrix}
        0 & 1 \\
        \frac12 & 0 & \frac12 \\
         & \frac12 & 0 & \ddots   \\
          & & \ddots & \ddots & \frac12 \\
         & &  & \frac12  & 0
    \end{pmatrix} - \frac{1}{2 b_K} 
    \begin{pmatrix}
        b_0 & \dots & b_{K-1}
    \end{pmatrix}
    \begin{pmatrix}
        0 \\ \vdots \\ 0 \\ 1
    \end{pmatrix}.
\end{equation*}
The eigenvalues $x_n$ of this matrix, found via \texttt{eigvals} in Julia, are the roots of the polynomial $q$.
In order to get roots of $p$, we simply solve the equation
\begin{equation*}
    \frac{z + z^{-1}}{2} = x_n
\end{equation*}
for both values of $z$. 
We note that this equation is sensitive near $x_n = \pm 1$, so the Chebyshev method here does not avoid larger errors of frequencies near $0$ and $1/2$.

\section{Choosing the Number of Circle Interpolation Modes}
\label{app:choosing-L}
Here, we give a method to determine the number of Fourier modes to project a given trajectory onto. 
We do this by controlling the condition number of the linear system \eqref{eq:circle-projection}.
We define the least-squares condition number of the linear system defined by the matrix $Y = W_T^{1/2} \Phi$ as
\begin{equation*}
    \kappa = \norm{Y^{\dagger}}_2 = \max_{\bm v} \frac{\norm{Y^{\dagger} \bm v}}{\norm{v}_2} = \frac{\sigma_1}{\sigma_L},
\end{equation*}
where $Y^{\dagger} = (\Phi^T W_T \Phi)^{-1}\Phi^T W_T^{1/2}$ is the Moore-Penrose pseudoinverse and $\sigma_1$ and $\sigma_L$ are the largest and smallest singular values of $Y$.

To bound the condition number, we first observe that the matrix $Y^T Y$ has a Toeplitz structure:
\begin{equation}
\label{eq:eta_n}
    (Y^T Y)_{mn} = \eta_{n-m} = \sum_{t=0}^{T}w_{t,T+1}\lambda_1^{(n-m)t},
\end{equation}
where we note that $\eta_0 = 1$. In this way, we have a Gershgorin circle bound on the eigenvalues $\sigma_j^2$ of $Y^T Y$ of
\begin{equation*}
    \abs{\sigma^2_j - 1} \leq 2 \gamma_L = 2 \sum_{n = 1}^{2L} \abs{\eta_n},
\end{equation*}
where we used the fact that $\abs{\eta_n} = \abs{\eta_{-n}}$. This can be translated to a bound on the condition number of
\begin{equation}
\label{eq:condition-bound}
    \kappa \leq \frac{\sqrt{1+2\gamma_L}}{\sqrt{1-2\gamma_L}}.
\end{equation}

The expression \eqref{eq:condition-bound} allows for an algorithm to bound the condition number of this linear system. We simply choose a maximum allowed radius $\gamma_{\max}$, and find the largest $L$ such that $\gamma_L < \gamma_{\max}$. This requires $\mathcal O(L T)$ operations, which is always cheaper than the linear solve. By default, we have found $\gamma_{\max} = 0.5$ to give good results. The algorithm is summarized in the following pseudocode:
\begin{algorithm}[H]
\caption{Parameterization Dimension $L$}
\label{alg:parameterization-dimension}
\begin{algorithmic}[1]
\Require System height $T$, Frequency $\omega$, Tolerance $\gamma_{\max}$
\State $L \gets \ceil{\frac{T-1}{2}}$
\State $\gamma_0 \gets 0$; $n \gets 0$
\While {$n \leq T$}
    \State $n \gets n + 1$
    \State $\gamma_n \gets \gamma_{n-1} + \eta_n$ via \eqref{eq:eta_n}
    \If {$\gamma_n \geq \gamma_{\max}$}
        \State \Return $L \gets \ceil{\frac{n-2}{2}}$
    \EndIf
\EndWhile
\Ensure Parameterization dimension $L$
\end{algorithmic}
\end{algorithm}
\section{Proof of Theorem \ref{thm:refpolyconvergence}}
\label{app:ideal-polynomial}
In this appendix, we show that for $d=1$ and a fixed regularity $M$, the RRE residual, and therefore the Birkhoff average, converges at a rate faster than the standard guarantee for the weighted Birkhoff average.
This improvement is primarily due to the use of the continued fraction representation of the rotation number rather than the Diophantine property. 
Using continued fractions, we find progressively better ``ideal polynomials,'' which exactly filter the leading frequencies of the signal.
Because these exactly match the frequencies (rather than the WBA filtering everything evenly), we approach the optimal rate given by the regularity.

The only restriction for the proof is that the continued fraction denominators are sufficiently close to each other, a fact which is given by a standard theorem. 
Using this, we then show that an ideal polynomial with a degree given by the continued fraction denominator is bounded. 
This bound gives enough information for the result.

\subsection{The distribution of continuous fraction denominators}
We begin by quoting a standard theorem on the distribution of continued fraction denominators. 
The continued fraction representation of a rotation number $\omega$ is given by the fraction
\begin{equation*}
    \omega = \frac{1}{c_1 + \frac{1}{c_2 + \ddots}}
\end{equation*}
where the coefficients $c_n \in \Zbb$ are positive integers.
By truncating the series for finitely many $c_n$, we obtain the convergents
\begin{equation*}
    \frac{N_1}{L_1} = \frac{1}{c_1}, \ \frac{N_2}{L_2} = \frac{1}{c_1 + \frac{1}{c_2}}, \ \frac{N_3}{L_3} = \frac{1}{c_1 + \frac{1}{c_2 + \frac{1}{c_3}}}, \ \dots
\end{equation*}
The continued fractions are best approximations of $\omega$, satisfying the inequality
\begin{equation}
\label{eq:bestapprox}
    \abs{\omega - \frac{N_n}{L_n}} \leq \frac{1}{L_n L_{n+1}} < \frac{1}{L_n^2}.
\end{equation}
This is the fundamental result which we will use for the proof of convergence.

Remarkably, one can prove that almost all continued fractions have denominators with similar limiting behavior: 
\begin{theorem}
\label{thm:denominator-convergence}
For almost all $\omega \in [0,1)$, the denominators of continued fractions $L_n$ obey the limit 
\begin{equation*}
    \lim_{n\to \infty}L_n^{1/n} \to \gamma,
\end{equation*}
where
\begin{equation*}
    \gamma = \exp\left(\frac{\pi^2}{12 \log 2}\right).
\end{equation*}
\end{theorem}
For an introduction to this theorem and more on the measure theory of continued fractions, a good reference is the book by Khinchin\cite{khinchin1997}.

We then transform this theorem to bounds on continued fraction denominators.
\begin{corollary}
\label{corollary:DenominatorDistribution}
    Let $r, s > 0$ and $\eta < 1$. For almost all $\omega \in [0,1)$, there is a value $M$ such that for all $m > M$, there exists a continued fraction approximation $p_n/L_n$ of $\omega$ such that 
    \begin{equation}
        \label{eq:denominatorbounds}
        r m^\eta < L_n < s m.
    \end{equation}
\end{corollary}
\begin{proof}
    First, we fix $\epsilon > 0$. Thm.~\ref{thm:denominator-convergence} tells us that there is an $N$ such that for all $n> N$
    \begin{equation*}
        \gamma - \epsilon < L_n^{1/n} < \gamma + \epsilon.
    \end{equation*}
    Combining the above equation with \eqref{eq:denominatorbounds}, it is clear that if we can show that for all $m>M$ we have
    \begin{align}
    \label{eq:bigdenominatorbound}
        r m^\eta < (\gamma - \epsilon)^n < L_n < (\gamma+\epsilon)^n < s m,
    \end{align}
    then we have proven our theorem.

    We note that $n>N$ must be satisfied if we choose $m$ large enough so that if the first inequality is satisfied, i.e.
    \begin{equation*}
        n > \frac{\log r + \eta \log m}{\log(\gamma - \epsilon)} > N.
    \end{equation*}
    Additionally, when
    \begin{equation*}
        \log s - \log r +(1- \eta) \log m > 2 \log \gamma > 2 \log (\gamma - \epsilon), 
    \end{equation*}
    there are at least two points $n$ that satisfy
    \begin{equation*}
        \log r + \eta \log m < n \log (\gamma - \epsilon) < \log s + \log m.
    \end{equation*}
    In particular, there is a point before the midpoint of the bound such that
    \begin{multline*}
        \log r + \eta \log m < n \log (\gamma - \epsilon) < \\ \frac{1}{2}(\log r + \log s) + \frac{1+ \eta}{2} \log m.
    \end{multline*}
    Using the above expression, we additionally know that
    \begin{multline*}
        n \log(\gamma + \epsilon) < \\ \frac{\log(\gamma + \epsilon)}{\log(\gamma - \epsilon)}\left[ \frac{1}{2}(\log r + \log s) + \frac{1+ \eta}{2} \log m \right].
    \end{multline*}
    If we choose $\epsilon$ small enough so that
    \begin{equation*}
        \frac{\log(\gamma + \epsilon)}{\log(\gamma - \epsilon)} = \frac{2 C}{1+\eta},
    \end{equation*}
    where $C < 1$, then the above inequality reduces to
    \begin{equation*}
        n \log(\gamma + \epsilon) <  \frac{C}{1+\eta}(\log r + \log s) + C\log m .
    \end{equation*}
    From here, it is clear that there exists an $M$ so that $m>M$ implies that 
    \begin{equation*}
        n \log(\gamma + \epsilon) < \log s + \log m,
    \end{equation*}
    giving our upper bound.
\end{proof}

\subsection{The Polynomial}
To prove Thm.~\ref{thm:refpolyconvergence}, we create a sequence of filters that annihilate a certain number of the frequencies $\lambda_n$. 
If we let $0 < \alpha < 1/4$ and $n \geq 1$, we do this by defining the reference polynomial for frequency $\omega$ and island period $p$ as
\begin{multline}
\label{eq:refpoly}
    \refpoly(z) = \prod_{j = 1}^{\floor{\alpha p L_{n}}}\frac{(z - \lambda_j)(z - \lambda_{-j})}{(1 - \lambda_j)(1 - \lambda_{-j})}\\
    \times \prod_{j = \floor{\alpha p L_{n}}+1}^{\floor{{p L_{n}}/{2}}}\frac{(z - \mu_j)(z - \mu_{-j})}{(1 - \mu_j)(1 - \mu_{-j})}
\end{multline}
where for $0 \leq \ell < p$ and $j\geq 0$
\begin{equation*}
    \mu_{j p + \ell} = e^{2 \pi i (j N_n/L_n + \ell)/p}.
\end{equation*}
A fraction $2\alpha$ of the roots of $\refpoly$ exactly match $\lambda_j$, and the rest of the roots $\mu_j$ lie on roots of unity. 
The polynomial has been chosen so that it is comparable to the polynomial with roots at roots of unity, which is well understood as a cardinal function of the discrete Fourier transform.

In particular, we have:
\begin{lemma}
\label{lemma:refpolybound}
    Let $0 < \omega \leq 1$ have a sequence of continued fractions approximants with $\{N_n / L_n \}$. For any $n \geq 1$ and $\alpha < 1/4$, the reference polynomial obeys the bound 
    \begin{equation*}
        \abs{\refpoly(z)} \leq C_\alpha
    \end{equation*} 
    for $\abs{z} \leq 1$ and $C_\alpha > 1$, where $C_\alpha$ is independent of $\omega$ and $n$.
\end{lemma}
\begin{proof}
    The proof consists of comparing the polynomial above to the polynomial with equispaced nodes 
    \begin{equation*}
        Q_{\alpha, n}(z) = \prod_{j = 1}^{\floor{{L_{n}}/{2}}}\frac{(z - \mu_j)(z - \mu_{-j})}{(1 - \mu_j)(1 - \mu_{-j})}.
    \end{equation*}
    We begin by assuming $L_n$ is odd. In this case, the polynomial coefficients of $Q_{\alpha, n}$ can be thought of as a column of the inverse discrete Fourier transform. This is because we have $Q(1) = 1$ and $Q(\mu_j) = 0$ for $1 \leq \abs{j} \leq \floor{\nicefrac{L_n}{2}}$ So, $Q_{\alpha, n}$ is simply the interpolation of a unit vector that is nonzero at $1$, i.e.~$\dbf = F^{-1} \ebf_0$ where $\dbf$ are the polynomial coefficients, $F_{mn} = \mu_{n-M_n}^m$ is the orthogonal discrete Fourier transform mapping from Fourier coefficients to values at equispaced nodes, and $M_n = \floor{p L_n/2}$.

    The problem of how Fourier interpolation changes when the nodes are perturbed is addressed in the paper by Yu and Townsend\cite{yu2023}. The authors prove an analogue to the Kadec-$1/4$ theorem \cite{kadets1964} from sampling theory, which states that as long as the locations of Fourier nodes does not move a length farther than $\alpha/4$ from the equispaced grid, the non-uniform Fourier transform matrix $\tilde{F}$ does not differ too much in norm from the equispace Fourier transform, i.e.
    \begin{equation*}
        \norm{\tilde F} \leq (1 + \phi_\alpha)\norm{F}_2,
    \end{equation*}
    where $\norm{F} = L_n$, $\phi_\alpha = 1 - \cos \pi \alpha + \sin \pi \alpha$, and
    \begin{equation*}
        \tilde{F}_{mn} = \begin{cases}
            \lambda_{n-M_n}, & \abs{n-M_n} \leq \floor{\alpha p L_{n}}, \\
            \mu_{n-M_n}^m, & \floor{\alpha p L_{n}} < \abs{n-M_n} \leq M_n.
        \end{cases}.
    \end{equation*} Additionally, via Weyl's inequality we have
    \begin{equation*}
        \norm{\tilde F^{-1}}_2 \leq \frac{1}{(1 - \phi_\alpha) \norm F}_2.
    \end{equation*}
    
    We can apply these identities to our problem. 
    Using the approximation bound \eqref{eq:bestapprox}, we find that $\abs{\lambda_j - \mu_j} < \alpha \pi / L_n$ for $\abs{j} < \floor{\alpha L_n}$, obeying the hypothesis for the $1/4$ theorem. 
    Hence, the coefficients of $\refpoly$ 
    \begin{equation*}
        \norm{\filter_{\alpha, n}}_2 = \norm{\tilde F^{-1} \ebf_0}_2 \leq \frac{1}{(1-\phi_\alpha)\sqrt{L_n}}. 
    \end{equation*}
    The polynomial evaluation bound directly follows from the above inequality. 
\end{proof}

\subsection{Convergence of the Least-Squares Residual}
With Lemma \ref{lemma:refpolybound} and Corollary \ref{corollary:DenominatorDistribution}, we have enough information to prove Theorem \ref{thm:refpolyconvergence}:
\begin{proof}[Proof of Theorem \ref{thm:refpolyconvergence}]
    Choose $r = s = 1/p$ and $\eta <1$. Then, from Corollary \ref{corollary:DenominatorDistribution} we know there is an $L$ such that $2K+1>L$ implies there is an $n$ such that the continued fraction representation of $\omega$ obeys
    \begin{equation*}
        (2K+1)^\eta < p L_n < 2K+1.
    \end{equation*}
    Therefore, the filter $\filter_{\alpha, n}$ with the coefficients of $\refpoly$ has a length less than the maximum allowed for the linear system. 
    To make the filter compatible with the full system, we pad $\filter_{\alpha, n}$ symmetrically with zeros on the top and bottom.
    We note that $\filter_{\alpha, n}$ obeys the palindromic and correct-mean constraints by construction.

    Now, consider the application of $\filter_{\alpha, n}$ on $U$. 
    Using the Fourier representation of $\dobsT$, we have
    \begin{equation*}
        U_{jk} = \dobsT((k+j)\omega) = \sum_{m \in \Zbb} \dobsT_m \lambda_m^{k+j}.
    \end{equation*}
    Applying the filter, we find
    \begin{align*}
        (U \filter_{\alpha, n})_j &= \sum_{m \in \Zbb} \dobsT_m \lambda_{m}^j \refpoly(\lambda_m),\\
            &= \sum_{\abs{m} > \floor{\alpha p L_{n}}} \dobsT_m \lambda_{m}^j \refpoly(\lambda_m).
    \end{align*}

    Using Lemma \ref{lemma:refpolybound}, we have
    \begin{align*}
        \abs{(U \filter_{\alpha, n})_j} &\leq C_\alpha \sum_{\abs{m} > \floor{\alpha p L_{n}}} \norm{\dobsT_m},\\
        &\leq C_\alpha C_M \left( \alpha p L_{n}\right)^{-M+1},\\
        &\leq C_\alpha C_M \alpha^{-M+1} (2K+1)^{\eta (-M+1)}.
    \end{align*}
    where we use the fact that $\dobsT \in C^M$ for the second inequality for some $C_M > 0$.
    Substituting into the full RRE quadratic form, we have
    \begin{align*}
        (U \filter_{\alpha,n})^T W_T (U \filter_{\alpha,n}) &= \sum_{t=0}^T w_{t,T}(U \filter_{\alpha, n})_j^2, \\
        & \leq C (2K+1)^{2\eta (-M+1)}.
    \end{align*}

    In the case that $\obsT$ is real-analytic, there exists a constant $\rho$ such that 
    \begin{equation*}
        \norm{\dobsT_m} < C_\rho \rho^{-m}.
    \end{equation*}
    This further implies that for some $\tilde{C} > 0$
    \begin{align*}
         \abs{(U \filter_{\alpha, n})_j} &< \tilde{C} \rho^{-\alpha p L_{n}}, \\
         &< \tilde{C} \rho^{-\alpha(2K+1)^{\eta}}
    \end{align*}
    After substitution into the quadratic form, we find our result.

    
\end{proof}

\bibliography{bibliography}

\end{document}